\documentclass[12pt]{article}
\newtheorem{theor}{Theorem}
\newtheorem{lemma}{Lemma}
\newtheorem{example}{Example}
\newtheorem{prop}{Proposition}

\newcommand{\aaa}{{\cal A}}
\newcommand{\bbb}{{\cal B}}

\newcommand{\ddd}{{\cal D}}

\newcommand{\ggg}{{\cal G}}

\newcommand{\lll}{{\cal L}}

\newcommand{\ppp}{{\cal P}}

\newcommand{\ttt}{{\cal T}}

\newcommand{\xxx}{{\cal X}}

\newcommand{\zzz}{{\cal Z}}

\newcommand{\xxb}{{\mbox{\boldmath $ x $}}}

\newcommand{\reff}[1]{(\ref{#1})}
\newcommand{\proof}{\vskip0.4cm \noindent {\bf Proof: }}
\newcommand{\eproof}{ \mbox{}\hfill$\sqcup\!\!\!\!\sqcap$ \vskip0.4cm\noindent}
 
\newcommand{\ie}{{\it i.e. }}

\newcommand{\virg}{\, , \, }

\newcommand{\ponto}{\, \cdot  \, }
 
\newcommand{\formu}[1]{\begin{eqnarray}\label{#1}}
\newcommand{\formub}{\begin{eqnarray} \nonumber}
\newcommand{\eformu}{\end{eqnarray}}

\newcommand{\re}{{I\!\!R}}

\newcommand{\setag}{\longrightarrow}

\newcommand{\pae}[2]{\stackrel{\hspace{-3mm}(\bf e)} {\Pi_{#1}^{#2}} }
\newcommand{\pam}[2]{\stackrel{\hspace{-3mm}(\bf m)} {\Pi_{#1}^{#2}} }

\newcommand{\It}{{\{ x : t \vert r_t (x) \vert > 1 \}}}
\newcommand{\Itc}{{\{ x : t \vert r_t (x)\vert \le 1 \}}}

\newcommand{\zeb}{\mbox{\boldmath $0$}}

\newcommand{\alphab}{\mbox{\boldmath $\alpha$}}

\begin{document}
\title{Inference Functions for Semiparametric Models}
\author{Rodrigo Labouriau 
 \thanks{Department of Mathematics,
               Aarhus University.
        }
      }
\date{Fall, 2020}
\maketitle

\begin{abstract}
The paper discusses inference techniques for semiparametric models based on suitable versions of inference functions. The text contains two parts. In the first part, we review the optimality theory for non-parametric models based on the notions of path differentiability and statistical functional differentiability. Those notions are adapted to the context of semiparametric models by applying the inference theory of statistical functionals to the functional that associates the value of the interest parameter to the corresponding probability measure. 
The second part of the paper discusses the theory of inference functions for semiparametric models. We define a class of regular inference functions, and provide two equivalent characterisations of those inference functions: One adapted from the classic theory of inference functions for parametric models,  and one motivated by differential geometric considerations concerning the statistical model. Those characterisations yield an optimality theory for estimation under semiparametric models. We present a necessary and sufficient condition for the coincidence of the bound for the concentration of estimators based on inference functions and the semiparametric Cram\`er-Rao bound. Projecting the score function for the parameter of interest on specially designed spaces of functions, we obtain optimal inference functions. Considering estimation when a sufficient statistic is present, we provide an alternative justification for the conditioning principle in a context of semiparametric models. The article closes with a characterisation of when the semiparametric Cram\`er-Rao bound is attained by estimators derived from regular inference functions.
\end{abstract}


\newpage

\tableofcontents
 
\newpage

\section{Introduction}

In this article, we revise the classic optimality theory for non- and semiparametric models.  A range of notions of path differentiability and tangent spaces are introduced and their inter-relations studied.
Next, it is studied some concepts of statistical functional differentiability.
Here the differentiability is considered relatively to a pointed cone contained in the tangent space and not relatively to the whole tangent
space, as is currently in the literature. These cones are referred to as the tangent cones. The optimality theory of differentiable functionals is
reviewed next. Again, the results are stated relative to the tangent cone and not with respect to the whole tangent space, as is usual.
The estimation of the interest parameter of semiparametric models is studied by applying the optimality theory to a specially designed functional called the interest parameter functional, which associates to any probability measure in the model in play the value of the interest parameter associated to it.
We will consider an increasing range of tangent cones. Here, the larger is the tangent cone used, the sharper is the bound for the concentration of regular estimators obtained. However, too large tangent cones may imply that the interest parameter functional is differentiable only under somehow stringent regularity conditions on the model.
We show how the imposition of such conditions usually done in the
literature can be avoided by using adequate choices of tangent cones.
The bound for the concentration for regular inferencesequences obtained with this choice of the tangent cone is referred to as the semiparametric Cram\'er-Rao bound.

\section{Path and Functional Differentiability} \label{cap2}

We consider in this section some aspects of the general theory of
non-parametric statistical models which will be useful for
the theory of semiparametric models.
The key notions introduced here are the path differentiability,
the associated concept of tangent spaces and tangent sets, and the notions of 
functional differentiability.

In section \ref{gents2} we study a range of concepts of path differentiability
and comparisons of those notions are provided.
An important point there is the equivalence between the
Hellinger differentiability, often used in the literature
(see Bickel {et al.}, 1993), and the weak differentiability
(see Pfanzagl, 1982, 1985 and 1990). Two auxiliary notions of path 
differentiability are introduced: strong and mean differentiability. 
It is proved that weak (or Hellinger differentiability) is an intermediate 
notion of path differentiability, weaker than strong differentiability and 
stronger than mean differentiability.
A new notion of path differentiability, called essential differentiability,
is introduced. We will interpret the 
tangents of  essential differentiable paths as score functions of 
one dimensional ``regular submodels'' in the classical sense. 
Since the essential differentiability
is weaker than the other notions provided, this interpretation extends
immediately  to all the other path differentiability notions considered.
 
In section \ref{gents3} some  differentiability notions of  functionals
are studied.
In the approach given a cone contained in the tangent set (\ie the class
of tangents of differentiable paths) is chosen and the differentiability 
of the functional in question will be defined relatively to this cone 
(termed tangent cone).
Alternative notions of functional differentiability are given by adopting
different notions of path differentiability and/or using different tangent 
cones. 
As we will see, the stronger the path differentiability notion used
and the smaller is the tangent cone,
the weaker the notion of differentiable functionals induced, in the sense
that more statistical functionals are differentiable.
We provide next some lower bounds for the concentration of ``regular''
sequences of estimators for a differentiable functional under a repeated
sampling scheme.
The weaker the path differentiability required and the larger
the tangent cone adopted, the sharper are the bounds obtained.
The theory will be applied to estimation in semiparametric models in section
\ref{ggg4}.

\subsection{Differentiable paths \label{gents2}}

The main purpose of this section is to introduce   
the mathematical machinery necessary to extend
the notion of score function, classically defined for 
parametric models, to a context where no (or only a partial)
finite dimensional parametric structure is assumed.
The key idea here is to consider one-dimensional
submodels of the  family $\ppp$ of probability measures
(typically infinite dimensional). 
These submodels will be called paths. 
Following the steps of Stein (1956), one should consider 
a class of submodels (or paths) sufficiently regular in 
order to have a score function well defined and well 
behaved for each submodel, in the sense that, at least,
each score function should be unbiased (\ie have 
expectation zero) and have finite variance.
Stein's idea is to use the worst possible regular submodel 
to assess the difficulty of statistical inference procedures
for the entire family $\ppp$.
Evidently, if the class of ''regular submodels'' is too small,
no sensible results are to be expected from that procedure.
On the other hand, if the class of ''regular submodels'' is
too large, the Stein's procedure can become intractable
or no simplification is really gained,
which is not in the spirit of the method proposed.
Hence, when applying the Stein procedure it is our task to 
find a class of ''regular submodels'' with the adequate size.

The idea of ''regular submodel'' mentioned will be formalised
by introducing the notion of path differentiability.
A range of concepts of path differentiability are studied
in this section, all of them fulfilling the
minimal requirement for a ''regular submodel'', \ie the 
score functions of the differentiable paths (viewed
as submodels) will be automatically well defined,
unbiased and possess finite variances.
The strongest notion of path differentiability
considered is the $L^\infty$ differentiability
(or pointwise differentiability) and the weakest
notion is the essential differentiability.
It will turn out that a notion of path differentiability called
``Hellinger differentiability''
 (or ``weak differentiability'') is the weakest notion that 
captures some important essential statistical properties of 
the model $\ppp$.
Another distinguished notion considered is the
$L^2$ differentiability which will involve calculations
with Hilbert spaces, simplifying all the computations required. 
The $L^2$ differentiability coincides with the
Hellinger differentiability in most of the examples considered in
this thesis.
It turns that the $L^2$ differentiability will be useful in the 
theory of inference functions.
 
This section is organised as follows.
Subsection \ref{gendiff1} studies the basic notion
of path differentiability and some general properties
of differentiable paths.
Some specific concepts of differentiability are introduced
in the  subsections  \ref{gendiff2}, \ref{gendiff3} and \ref{gendiff4} 
where weak or Hellinger, $L^q$ and essential differentiability are studied, 
respectively.
The associated notions of tangent sets and tangent spaces
are discussed in subsection \ref{gendiff5}.

\subsubsection{General definition of path differentiability}
\label{gendiff1}

We give next a more precise definition of the terms ''submodel''
and ''regular submodel'' informally used in the previous discussion.
Recall that we were interested in defining a one-dimensional
submodel contained in the family $\ppp$ for which the score function
would be well defined and well behaved.

Let us consider a subset $V$ of $[0,\infty )$ which contains zero and
for which zero is an accumulation point.
The set $V$ will play the role of the parameter space in the 
''submodel'' we define.
Typical examples are: 
$[0,\epsilon )$ for some $\epsilon >0$
and $\{ 1/n : n\in N \}\cup\{ 0 \}$.
A mapping from $V$ into $\ppp^*$ assuming the value $p\in\ppp^*$
at zero is said to be a {\it path \/} converging to $p$.
Here the image of $V$ under a path plays the role of the 
''submodel'' of $\ppp$ and the path acts as a one-dimensional 
parametrisation of the ''submodel''.
It is convenient to represent a path by a generalised sequence
$\{ p_t\}_{t\in V} = \{ p_t\}$, where for each $t\in V$, 
$p_t\in\ppp^*$ is the value of the path at $t$.

We introduce next the notion of differentiability which will enable us
to formalise more precisely what in the Stein program is the class
of ''{\it regular\/} submodels''.  
A path $\{ p_t\}_{t\in V}$ (converging to $p$)
is {\it differentiable} at $p\in\ppp^*$
if for each $t\in V$ we have the representation
\formu{gendiff101}
p_t (\ponto ) = p (\ponto ) + t p (\ponto ) \nu (\ponto )
                            + t p (\ponto ) r_t  (\ponto )
\eformu
for a certain $\nu (\ponto ) \in L^2_0 (p)$, and
\formu{gendiff102}
r_t \longrightarrow 0 , \,\,\, \mbox{ as } t\downarrow 0 
\, .
\eformu
The convergence in (\ref{gendiff102}) is
in some appropriate sense to be specified later.
In fact, in the next subsections we explore several notions of
path differentiability by introducing alternative 
definitions for that convergence. 
The term $r_t$ in (\ref{gendiff101}) will be referred to as
the {\it remainder term\/}. 

The function $\nu :\xxx\longrightarrow\re$ given in
(\ref{gendiff101}) is said to be the {\it tangent} associated
to the differentiable path $\{ p_t\}$.
Here the tangent plays the role of the score function of the 
submodel parametrised by $t\in V$ at $p_0=p$.
To see the analogy with the score function suppose that the
convergence of $r_t$ in (\ref{gendiff102}) is in the sense of 
the pointwise convergence. In that case the tangent 
coincides with the score function of the submodel associated
with the differentiable path $\{ p_t\}$ at $p_0=p$.
In the general case, where the convergence of $r_t$ is not 
necessarily pointwise convergence, the general chain rule
for differentiation of functions in metric spaces
(see Dieudonn\' e , 1960)
can often be applied to justify our interpretation of the 
tangent.  
We stress that according to our definition, the tangent of a 
differentiable path (or alternatively the score of a regular
submodel) has automatically finite variance and mean zero
(\ie it is in $L^2_0 (p)$). 

Before embracing the study of  notions of differentiability
generated by some specific definitions of the convergence of $r_t$,
we give a useful and trivial general property of 
remainder terms of differentiable paths.
Suppose that a path $\{ p_t\}$ is differentiable 
at $p\in\ppp^*$ with representation given by (\ref{gendiff101}), with 
$\nu\in L^2_0 (p)$.
Then we have,  for each $t\in V$
\formu{gendif17a}
r_t(\ponto ) =
 \frac{p_t (\ponto )- p(\ponto )}{t p(\ponto )} - \nu (\ponto )
\eformu
and
\formu{gendif116a}
\nonumber
\int_\xxx r_t (x ) p(x) \lambda (dx) = 
\int_\xxx
\left\{ 
 \frac{p_t (x )- p(x )}{t p(x )} - \nu (x )
\right \}
p(x) \lambda (dx) = 0
\, .
\eformu

\subsubsection{Hellinger and weak path differentiability}
\label{gendiff2}

Most of the estimation theory for non- and semi-parametric models found in
the literature (see Bickel {\it et al.\/}, 1993 and references therein) is
developed using the notion of Hellinger differentiability studied next.
This notion appears in the literature in two equivalent forms:
weak differentiability (see Pfanzagl 1982, 1985 and 1990) and Hellinger
differentiability (see H\' ajeck, 1962, LeCam, 1966 and Bickel {\it et al.\/},
1993).
This notion of differentiability plays a central role in the 
theory presented because it enables us to grasp some essential statistical
properties of the models considered.
For instance, the Hellinger differentiability is equivalent to local 
asymptotic normality 
of the submodel defined by the path.
Moreover, the Hellinger differentiability is used in the so called convolution
theorem, which gives a bound for the concentration of a rich class
of estimators (the regular asymptotic linear estimators).

We begin by introducing the weak differentiability which is in the
general form of path differentiability formulated before.
A path $\{ p_t\}_{t\in V}$ is {\it weakly differentiable \/} at $p\in\ppp^*$
if  there exist $\nu\in L^2_0 (p)$ and a generalised sequence of functions
$\{ r_t \}_{t\in V}$ such that for each $t\in V$
\formu{gendiff19}
\nonumber
p_t (\ponto ) = p (\ponto ) + t p (\ponto ) \nu (\ponto )
                            + t p (\ponto ) r_t  (\ponto )
\eformu
and
\formu{weak1}
\frac{1}{t}
\int_\It 
\vert r_t (x)  \vert p(x) \lambda (dx) \longrightarrow 0
\, , \,\, \mbox{ as } t\downarrow 0
\, ,
\eformu
\formu{weak2}
\int_\Itc 
\vert r_t (x)  \vert^2 p(x) \lambda (dx) \longrightarrow 0
\, , \,\, \mbox{ as } t\downarrow 0
\, .
\eformu
In other words, $\{ p_t\}$ is weakly differentiable if it is differentiable
according to the general definition of path differentiability with the
convergence of the generalised sequence $\{ r_t\}$ given by 
(\ref{weak1}) and (\ref{weak2}).

Let us introduce now the Hellinger differentiability of paths.
The key idea in this approach is to characterise the family $\ppp$
of probability measures by the class of square roots of the densities,
instead of the densities.
The advantage of this alternative characterisation is that the square
roots of the densities are in the Hilbert space
\formu{gendiff20}
\nonumber
L^2(\lambda ) = \left \{ 
                    f :\xxx\longrightarrow\re \, : 
               \int_\xxx f^2(x) \lambda (dx) < \infty 
            \right \}
\, .
\eformu
In this way the statistical model in play is naturally embedded into a
space with a rich mathematical structure.
Using the usual topology of $L^2(\lambda )$ one defines the differentiability of 
paths in the sense of Fr\' echet (or in this case, since the domain of the path
is contained in $\re$, the equivalent notions of Hadamard and Gateaux
differentiability could be used also).
The precise definition of Hellinger differentiability is the following.
A path  $\{ p_t\}_{t\in V}$ 
is {\it Hellinger differentiable} at $p\in\ppp^*$ if 
there exists a generalised sequence  $\{ s_t\}_{t\in V}$ 
in $L^2_0 (p)$ converging to zero as $t\downarrow 0$, \ie
\formu{gendiff31}
\| s_t \|_p \longrightarrow 0 \, , \, \mbox{ as } t \downarrow 0
\eformu
and $\nu\in L^2_0(p)$ such that
\formu{gendiff32}
p_t ^{1/2}(\ponto ) = p ^{1/2} (\ponto ) 
                      + t p  ^{1/2}(\ponto ) \frac{1}{2} \nu (\ponto )
                      + t p ^{1/2} (\ponto )  s_t  (\ponto )
\,  .
\eformu
The factor $\frac{1}{2}$ in the second term of the right  side
of (\ref{gendiff32}) will serve to accommodate with the other notions
of differentiability.
Note that each $s_t$ is in fact in $L^2_0 (p)$.
For, from (\ref{gendiff32})
\formu{gendiff33}
s_t (\ponto) = \frac{p_t^{1/2} (\ponto ) - p^{1/2} (\ponto )}
                    {t p^{1/2} (\ponto )}
             - \frac{\nu (\ponto )}{2} 
\, .
\eformu
Since 
$
\int_\xxx \left \{ \frac{p_t^{1/2} (x)}{p^{1/2} (x)} \right \}^2
          p(x) \lambda (dx)
=
\int_\xxx p_t(x) \lambda (dx) 
= 1 < \infty 
\, ,
$
we have that $p_t^{1/2} (\ponto ) / p^{1/2} (\ponto )\in L^2 (p)$,
and hence 
$
\frac{p_t^{1/2} (x)}{t p^{1/2} (x)}
= \frac{1}{t} \left \{ \frac{p_t^{1/2} (x)}{p^{1/2} (x)} - 1 \right \}
\in L^2 (p)
$.

\begin{prop}
\label{gendiffp0}
A path $\{ p_t\}$ is Hellinger differentiable 
 if and only if $\{ p_t\}$ is weak differentiable.
\end{prop}
\proof
See Pfanzagl (1985).
\eproof

\subsubsection{$L^q$ path differentiability}\label{gendiff3}

We study next a useful range of path differentiability notions.
These notions will serve us to graduate how strong is the Hellinger or 
weak differentiability; and they will be used auxiliary in the 
calculation of the weak tangents of weak differentiable paths.
In spite of the secondary role these differentiability notions
play in our development, they are important in the general theory
of differentiability of statistical functionals, in particular
in the theory of von Mises functionals.
The $L^2$ differentiability defined below will be useful when studying
the use of inference functions for semiparametric models.

The main idea here is to consider the $L^q$ convergence for the 
generalised sequence $\{ r_t \}$ appearing in the definition of
differentiable paths. The precise definition is the following.
A path  $\{ p_t\}_{t\in V}\subseteq \ppp^*$ is 
{\it $L^q$ differentiable\/} at $p\in\ppp^*$, for $q\in [1,\infty]$,
if there exist $\nu\in L^2_0(p)$ and a generalised sequence
$\{ r_t \}$ in $L^q(p)$ such that for each $t\in V$,
\formu{gendiff40}
p_t (\ponto ) = p (\ponto ) + t p (\ponto ) \nu (\ponto )
                            + t p (\ponto ) r_t  (\ponto )
\eformu
and
\formu{gendiff41}
\| r_t \|_{L^q(p)} \longrightarrow 0 \, , \,\, \mbox{ as } t\downarrow 0 
\, .
\eformu
The following proposition relates the notions of $L^q$ path differentiability.
\begin{prop}
\label{gendiffp1}
Consider $r,q\in [1,\infty]$ such that $r\le q$.
If a path is $L^q$ differentiable at $p\in\ppp^*$, then it is 
also $L^r$ differentiable at $p$ with the same tangent.
\end{prop}
\proof
The proposition follows immediately from the fact that 
convergence in $L^q(p)$ implies convergence in $L^r(p)$.
\eproof
 
There are two distinguished cases of $L^q$ path differentiability:
strong and mean differentiability corresponding to $L^2$ and $L^1$
differentiability respectively.
The $L^1$ differentiability is remarkable because it is the weakest 
notion of differentiability found in the literature, and the $L^2$
differentiability distinguish itself because the $L^2$ spaces,
when endowed with the natural inner product, are Hilbert spaces,
which simplifies significantly the calculations.

We study next the relation between weak and $L^q$ path differentiability.
As we will see in the  propositions \ref{gendiffp2} and \ref{gendiffp3}
given above, weak differentiability is an
intermediate notion of path differentiability between $L^2$ and $L^1$ differentiability.
\begin{prop}
\label{gendiffp2}
If a path is $L^2$ differentiable at $p\in\ppp^*$,
then it is weakly (or Hellinger) differentiable at $p$,
with the same tangent.
\end{prop}
\proof
Let $\{ p_t\}$ be a differentiable path in the $L^2$ sense
with representation (\ref{gendiff40})
and $\| r_t \|_{L^2(p)} \longrightarrow 0$ as $t\downarrow 0$.
We show that the path $\{ p_t\}$ fulfills the conditions
(\ref{weak1}) and (\ref{weak2}) for the convergence of the remainder
term in the sense of the weak path differentiability.
For,
\formu{genbb12}
\nonumber
\hspace{-10mm}
\frac{1}{t} \int_\It 
 \vert r_t (x) \vert p(x) \lambda (dx)
& \le &
\frac{1}{t} \int_\It 
 t \vert r_t (x) \vert  \vert r_t (x) \vert p(x) \lambda (dx)
\\ \nonumber
& = & 
 \int_\It 
  \vert r_t (x) \vert^2  p(x) \lambda (dx)
\\ \nonumber
& \le &
\int_\xxx \vert r_t (x) \vert^2  p(x) \lambda (dx)
\\ \nonumber
& = &
\|  r_t \|_p^2 \longrightarrow 0 
\, , \,\, \mbox{ as } t\downarrow 0
\, .
\eformu
Hence $\{ r_t\}$ satisfies (\ref{weak1}).
On the other hand,
\formu{genbb11}
\nonumber
\hspace{-9.5 mm}
\int_\Itc \vert r_t (x) \vert^2 p(x) \lambda (dx)
\le
\int_\xxx \vert r_t (x) \vert^2 p(x) \lambda (dx)
=
\|  r_t \|_p^2 \longrightarrow 0 
 , \mbox{ as } t\downarrow 0
 .
\eformu
Hence $\{r_t \}$ satisfies (\ref{weak2}). We conclude 
that $\{ p_t\}$ is differentiable in the weak sense with tangent $\nu$.
\eproof
\begin{prop}
\label{gendiffp3}
If a path is weakly (or Hellinger) differentiable at $p\in\ppp^*$,
then it is $L^1$ differentiable (or differentiable in mean) at $p$,
with the same tangent.
\end{prop}
\proof
Take a path $\{p_t \}$ weakly differentiable at $p$ with tangent
$\nu$. There exists a generalised sequence of functions $\{ r_t \}$
satisfying (\ref{weak1}) and  (\ref{weak2}),
such that for all $t\in V$,
\formu{gendiff50}
\nonumber
p_t (\ponto ) = p (\ponto ) + t p (\ponto ) \nu (\ponto )
                            + t p (\ponto ) r_t  (\ponto )
\, .
\eformu
Note that  (\ref{weak1}) implies that 
\formu{gendiff51}
\int_\It \vert r_t (x) \vert p(x) \lambda (dx)
\longrightarrow 0 
\, , \mbox{ as } t \downarrow 0
\, ,
\eformu
and  (\ref{weak2}) implies that 
\formu{gendiff52}
\int_\Itc \vert r_t (x) \vert p(x) \lambda (dx)
\longrightarrow 0 
\, , \mbox{ as } t \downarrow 0
\, .
\eformu
For, (\ref{weak2}) is equivalent to $L^2(p)$ convergence of
$s_t (\ponto ) := r_t (\ponto ) \chi_\Itc (\ponto)$ to zero.
From (\ref{gendiff41}),
\formub
\int_\Itc \vert r_t (x) \vert p(x) \lambda (dx)
= \| s_t \|_{L^1(p)} \le \| s_t \|_{L^2(p)}
\longrightarrow 0
\, .
\eformu

Combining (\ref{gendiff51}) and (\ref{gendiff52}) we obtain
\formu{gendiff53}
\nonumber
\int_\xxx \vert r_t (x) \vert p(x) \lambda (dx)
& = &
\int_\Itc \vert r_t (x) \vert p(x) \lambda (dx)
\\ \nonumber & & 
+
\int_\It \vert r_t (x) \vert p(x) \lambda (dx)
\longrightarrow 0 
\, , \mbox{ as } t \downarrow 0
\, .
\eformu
We conclude that $\{ r_t \}$ is $L^1$ differentiable at $p$ with 
tangent $\nu$.
\eproof

\subsubsection{Essential path differentiability }
\label{gendiff4}

We study next the weakest notion of path differentiability considered
in this text.
A path $\{ p_t \}_{t\in V}\subseteq \ppp^*$ is {\it essential differentiable \/}
at $p\in\ppp^*$ if there exists $\nu\in L^2_0 (p)$ and a generalised sequence
$\{ r_t :\xxx\longrightarrow \re \}_{t\in V}$ of $(\aaa , \bbb (\re ) )$-
measurable functions such that for each $t\in V$,
\formu{gendiff701}
p_t (\ponto ) 
=
p (\ponto ) + t p(\ponto )\nu (\ponto ) 
            + t p(\ponto )r_t (\ponto )
\eformu
and for any sequence $\{ k_n \}_{n\in N}\subseteq V $
such that $k_n\longrightarrow 0$ as 
$n\rightarrow\infty$ there is a subsequence 
$\{ k_i \}_{i\in N}\subseteq \{ k_n \}_{n\in N}$
such that $r_{k_i} (\ponto ) \longrightarrow 0$
$p$-almost surely as $i\rightarrow\infty$.

We show next that essential differentiability is weaker than differentiability
in mean which, in view of propositions \ref{gendiffp2} and \ref{gendiffp3}
implies that the essential differentiability is the weakest notion of path
differentiability considered here.
\begin{prop}
\label{gendiffp5}
If a path is $L^1$ differentiable, then it is essential differentiable,
with the same tangent.
\end{prop}
\proof
The generalised sequence $\{ r_t \}_{t\in V}$ is Cauchy,
because it converges in $L^1$ to zero.
Using theorem 3.12 in Rudin (1987, page 68)
the essential differentiability follows.
\eproof

\newpage

The following scheme represents the interrelation between the various
notions of path differentiability considered.
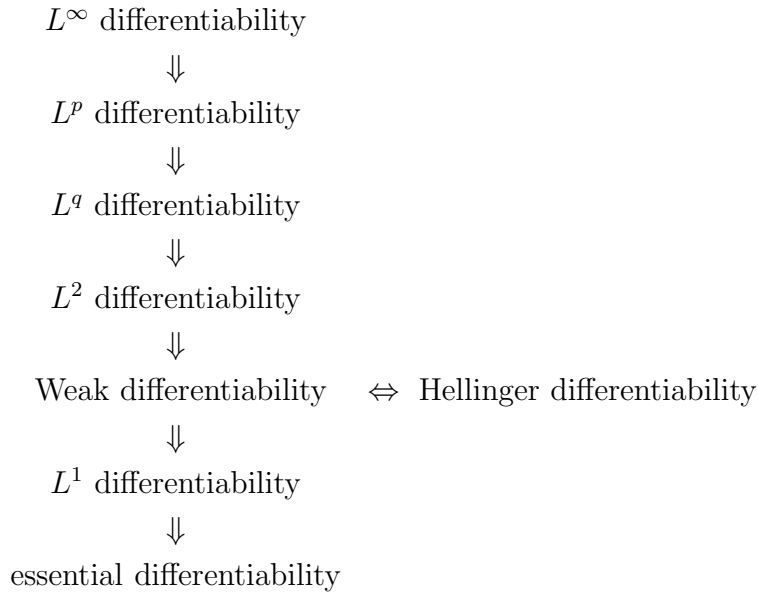
\formub
 & L^\infty \mbox{ differentiability} &
\\ \nonumber
& \Downarrow &
\\ \nonumber
 & L^p \mbox{ differentiability} &
\\ \nonumber & \Downarrow &
\\ \nonumber
 & L^q \mbox{ differentiability} &
\\ \nonumber &  \Downarrow &
\\ \nonumber
 & L^2 \mbox{ differentiability} &
\\ \nonumber & \Downarrow &
\\ \nonumber
 & \mbox{ Weak differentiability} & \Leftrightarrow
\mbox{ Hellinger differentiability}
\\ \nonumber  & \Downarrow &
\\ \nonumber
 & L^1 \mbox{ differentiability} &
\\ \nonumber  & \Downarrow &
\\ \nonumber
 & \mbox{essential differentiability} &
\eformu
Here $2<p<r<\infty$.

\subsubsection{Tangent spaces and tangent sets}
\label{gendiff5}

Re-taking the Stein approach, the notion of differentiable path formalised
the idea of ''regular one-dimensional submodel'', the tangent of a 
differentiable path playing the role of the score function of these submodels.
Here we elaborate the notion of tangent set which is the class of all possible
tangents of differentiable paths.
This will be useful to work with the idea of ''worst possible case''
contained informally in the Stein method, and to specify global
properties common to all the scores of ''regular one-dimensional
submodels''.
For technical reasons we need in fact to work in many situations
with the smallest closed subspace containing the tangent set, which
is called the tangent space.

In the next section we will define a notion of differentiability 
for statistical functionals. There the tangent set will play the role
of ''test functions'', analogous to the role of test functions when one 
defines the differentiability of tempered distributions
(see Rudin, 1973).
The notion of tangent space plays a crucial role when studying the 
theory of models with nuisance parameters. 
There we will need to obtain a component of a partial score
function orthogonal (in the $L^2$ sense, \ie uncorrelated) to the 
scores of a model obtained by fixing the parameter of interest and
letting the nuisance parameter vary. This component of the partial score 
function is obtained by orthogonal projection of the score function onto
the orthogonal complement of the tangent space (or nuisance tangent 
space as we will call the tangent space of the submodel we mentioned).
It will then be comfortable to work with a closed subspace of $L^2$.
We remark that it can be proved that the tangent set is a pointed cone,
but in general not even a vector space. Therefore the necessity to 
introduce the notion of tangent space as given here.

The formal definition of tangent space and tangent set depends on the 
notion of path differentiability one uses.
We give next a general definition of tangent set and tangent space 
which will be made precise when we specify the notion of path
differentiability we use.
Suppose we adopt a certain definition of path differentiability
according to which a differentiable path at $p\in\ppp^*$, 
say $\{ p_t \}$, has representation, for each $t\in V$,
\formu{gendiff801}
p_t (\ponto ) 
=
p (\ponto ) + t p(\ponto )\nu (\ponto ) 
            + t p(\ponto )r_t (\ponto )
\eformu
and 
\formu{gendiff802}
r_t \longrightarrow 0 
\, , \mbox{ as } t\downarrow 0
\, ,
\eformu
where the convergence in (\ref{gendiff802}) is in a certain sense known.
Then the {\it tangent set\/} of $\ppp$ at $p\in\ppp^*$ is the class
\formu{gendiff803}
\nonumber
\hspace{-10mm}
T^o (p) = T^o (p, \ppp ) =
\left \{
\begin{array}{r}
\nu \in L^2_0 (p )  : \, \exists V , \{p_t \}_{t\in V} \subseteq \ppp^* ,
                                       \{r_t \}_{t\in V} ,
\mbox{ such that } 
\\
\forall t\in V ,
\mbox{ (\ref{gendiff801}) and (\ref{gendiff802}) hold}
\end{array}
\right \} 
.
\eformu
The {\it tangent space\/} of $\ppp$ at $p\in\ppp^*$ is given by
\formu{gendiff803b}
\nonumber
T(p) = T(p,\ppp ) =
cl_{L^2_0 (p)} [ span \{ T^o (p, \ppp ) \} ]
\, .
\eformu
Since the tangent sets and spaces depend on the notion of path 
differentiability adopted, we speak of $L^q$ (for $q\in [1,\infty ]$), 
weak (or Hellinger) tangent sets and tangent spaces.
When necessary we use the notation  
$\stackrel{W} T$ for the weak   tangent space.
The $L^q$ tangent spaces are represented by
$\stackrel{q} T$ and the essential tangent spaces by 
$\stackrel{e} T$.

The following proposition relates the notions of tangent sets and
tangent spaces given.
\begin{prop}
\label{gendiffp7}
For each $p\in \ppp^*$ and for $2<q<r<\infty$ we have:
\formub
\stackrel{\infty} T (p) \subseteq 
\stackrel{r} T (p) \subseteq 
\stackrel{q} T (p) \subseteq
\stackrel{2} T (p) \subseteq
\stackrel{W} T (p) \subseteq
\stackrel{1} T (p) \subseteq 
\stackrel{e} T (p)
\, .
\eformu  
\end{prop}
\proof
Straightforward from the interrelations between the
notions of path differentiability.
\eproof

We close this section with two examples of the calculation of
tangent spaces.
\begin{example}

{\bf (Full tangent spaces of a large class of distributions)}
Consider the class $\ppp$ of all distributions in $\re$
dominated by the Lebesgue measure with continuous density
(with respect to the Lebesgue measure) and with support
(of the density) equal to the whole real line.
Denote the class of densities of $\ppp$ by $\ppp^*$.
We calculate the tangent space of $\ppp$ at each $p\in\ppp^*$.

Take an arbitrary element $\nu$ of $C_b\cap L^2_0(p)$.
Here $C_b$ denotes the class of continuous compact supported 
functions from $\re$ to $\re$.
It is a classical result of analysis that $C_b$ is dense
in $L^2 (p)$ (see Rudin, 1966), hence $C_b\cap L^2_0(p)$ is dense
in $L^2_0(p)$.
We show that $\nu\in T^0 (p,\ppp)$
(for {\it any\/} notion of tangent sets defined before).
Consider the path $\{ p_t \}$ given for $t\in [0,\infty)$
small enough, by
\formu{tangestr}
p_t(\ponto ) = p (\ponto ) + t p(\ponto )\nu(\ponto )
\, .
\eformu
We claim that for $t$ sufficiently small, $p_t\in\ppp^*$, which implies
that $\nu\in T^0 (p,\ppp )$.
It suffices to verify that $p_t$ is positive and integrates $1$.
For $t$ small $p_t$ is positive because $\nu$ is bounded and $p$ is
bounded in the support of $\nu$, hence the second term in the right hand of 
(\ref{tangestr}) is smaller than $p$ (for $t$ small).
That $p_t$ integrates $1$ follows from the fact that $\nu$ has 
expectation zero (with respect to $p$).
\eproof
\end{example}

It is not surprising that the previous enormous class of distributions
possesses a ``full'' tangent space.
The next example show that this could be the case even in families
where we have a lot of information about the distributions of
the family.
\begin{example}

{\bf (Full tangent space for families with information
on the moments)}
Consider the class $\ppp$ of all distributions in $\re$
dominated by the Lebesgue measure with continuous density
(with respect to the Lebesgue measure) and with support
(of the density) equal to the whole real line. Suppose further that 
the moments of all orders exist and that there exist a $\delta>0$,
a $k\in N$ and the constants $m_1,\dots ,m_k$ such that for each $i\in 
\{ 1, \dots , k \}$ the moment of order $i$ is contained in the open
interval $(m_i -\delta , m_i +\delta )$.
I claim that the tangent space of $\ppp$ at any $\ppp^*$ is 
$L^2_0 (p)$.
The proof follows the same line of the argument as given in the previous 
example.
Take a path as in (\ref{tangestr}) with $\nu\in C_b\cap L^2_0 (p)$.
For $t$ sufficiently small, $p_t$ will be positive, integrate to one, possess
finite moments of all orders, and the moments of order $i$, for $i\le k$
will be contained in the interval  $(m_i -\delta , m_i +\delta )$.
\eproof
\end{example}


\newpage
\subsection{Functional differentiability} 
\label{gents3}     

\subsubsection{Definition and first properties of functional differentiability}

We consider in this section a functional $\phi :\ppp^*\longrightarrow \re^q$
(for some $q\in N$)
which will play the role of a parameter of interest that we want to estimate.
Typical examples are the mean and the second moment functionals defined by 
$\phi(p) = \int_\xxx x p(x) \lambda (dx)$ 
and $\phi(p) = \int_\xxx x^2 p(x) \lambda (dx)$ respectively.
An important non trivial example for the theory of semiparametric models is 
the interest parameter functional defined next and studied in detail in 
section \ref{ggg4}.
\begin{example}
\label{semi1e}

{\bf Semiparametric models}

\noindent
Suppose that the family $\ppp^*$ of probability densities with 
respect to a measure $\lambda$ can be represented in the form
\formub
\ppp^* = \left  \{ p(\ponto ;\theta , z) \, : \, \theta\in\Theta
                  \subseteq\re^q
                  , z\in \zzz
         \right \}
\, .
\eformu
Here it is assumed that the mapping $(\theta , z)\mapsto p(\ponto ;\theta , z)$
is a bijection between $\Theta\times\zzz$ and $\ppp^*$.
The {\em interest parameter functional } $\phi :\ppp^*\longrightarrow \re^q$ 
is defined, for each $p(\ponto ;\theta , z)\in\ppp^*$, by
\formub
\phi\{ p(\ponto ; \theta , z) \} = \theta 
\, .
\eformu
\eproof
\end{example}

We introduce next a notion of functional differentiability that will enable
us to develop a theory of estimation for the functional $\phi$.
Let $p$ be a fixed element of $\ppp^*$. Consider a non empty subset 
$\ttt (p)$ of the tangent space at $p$.
A functional $\phi : \ppp^* \setag \re^q$ is said to be 
{\it differentiable at}
$p\in\ppp^*$ with respect to $\ttt (p)$ if there exists a function 
$\phi^\bullet_p :\xxx\setag\re^q$,
such that $\phi^\bullet_p \in \{ L^2_0 (p)\}^q$ 
and for each $\nu\in\ttt (p)$ there is a differentiable path $\{ p_t \}$ 
with tangent $\nu$ and 
\formu{bgr1}
\frac{\phi(p_t) - \phi (p)}{t}
\quad \setag \quad
<\phi^\bullet_p , \nu >_p \; , \quad \mbox{as } t\downarrow 0
\, .
\eformu
Here $<\phi^\bullet_p , \nu >_p$ is the vector with components given by
the inner product of the $q$ components of $\phi^\bullet_p$ and $\nu$.
The function $\phi^\bullet_p :\xxx\setag\re$ is said to be a {\it gradient}
of the functional $\phi$ at $p$ (with respect to $\ttt (p)$).
Note that $\phi^\bullet_p$ depends on the point $p$ at which we study the
differentiability of the functional $\phi$.
If a functional $\phi$ is differentiable at each $p\in\ppp^*$ we say that
$\phi$ is {\it differentiable}.

Since the definition of functional differentiability depends on the notion
of path differentiability, we speak of $L^\infty$, $L^p$, strong ($L^2$), 
weak, mean ($L^1$) and essential functional differentiability.
When necessary we superpose a symbol indicating the notion of path differentiability
in play.
When we are speaking generically or when it is clear from the context which 
notion of path differentiability is in play, we just use the notation
$\phi^\bullet_p$ for the gradient and $T^0 (p , \ppp^*)=T^0 (p)$,  
$T(p , \ppp^*)=T(p)$
for  the tangent set and the tangent space of $\ppp^*$ at $p$ respectively.

Note that the notion of functional differentiability introduced here
involves a subset $\ttt (p)$  of the tangent space and not 
necessarily the whole tangent space as is current in the literature.
This will give much more flexibility to the estimation theory developed.
Clearly the smaller is the class $\ttt (P)$ (or the stronger is the notion 
of path differentiability)
used, the weaker is the related functional differentiability.
On the other hand, the larger is the class $\ttt (p)$, the sharper
will be the results of the estimation theory related, in the sense
that the bounds for the lower asymptotic variance will be larger
or the optimality results will include more estimating sequences.
In this sense the ideal would be to choose the larger $\ttt (p)$
(and the stronger path differentiability) that makes differentiable
the functional under study.
Of course, we will have to require some mathematical properties for the classes
$\ttt (p)$ in order to obtain a notion of functional differentiability
useful for the estimation theory of differentiable functionals.
For instance, it will be assumed through (and silently) that 
$\ttt (p ) $ is a pointed cone (\ie if $\nu\in\ttt (p)$, then 
for each $\alpha\in \re_+ \cup \{ 0 \}$, $\alpha\nu\in\ttt (p)$).
We will refer form now on to $\ttt (p)$ as the {\it tangent cone}.
It will be necessary sometimes to require the tangent cones to be
convex.

We consider next a trivial example that illustrates the mechanics of the 
functional differentiability.
\begin{example}[Mean functional]
\label{exgr2}
Let $\lambda$ be a $\sigma$-finite measure defined on a measurable 
space $(\xxx ,\aaa )$.
Consider a family of probability measures $\ppp $ on 
$(\xxx ,\aaa )$ dominated by $\lambda$ given by the following representation
\formu{grzuca}
\ppp
= 
\left \{ 
\frac{dP}{d\lambda } (\ponto ) = p (\ponto ) : \; (\ref{bgr3})-(\ref{bgr6}) 
\mbox{ hold } 
\right \}
\, .
\eformu
The conditions to define $\ppp$ are
\formu{bgr3}
\forall x \in \xxx \, , \quad p(x)>0 \, ;
\eformu
\formu{bgr4}
\int_\xxx p(x) \lambda (dx) = 1 \, ;
\eformu
\formu{bgr5}
p \mbox{ is continuous} \, ;
\eformu
\formu{bgr6}
\int_\xxx x^2 p(x) \lambda (dx) \, \in \, \re_+ \, .
\eformu
We denote the class of densities of the elements of $\ppp$ with 
respect to $\lambda$ by $\ppp^* $.
Define the functional $M :\ppp^*  \setag \re$ by, for each 
$p\in\ppp^*$
\formub
M(p) = \int_\xxx x \, p(x) \lambda (dx)
\, .
\eformu 
We prove that $M$ is a differentiable functional with respect to the 
$L^2$ tangent space.
As we have seen in the previous section the tangent space of $\ppp$ at any
$p\in\ppp^*$ is the whole space $L^2_0 (p)$.

Take $p\in\ppp^* $ fixed and an arbitrary $L^2$-differentiable
path at $p$, say $\{ p_t \}_{t\in V}$, with representation given by
for each $t\in V$
\formub
p_t (\ponto ) = p(\ponto ) + t p (\ponto ) \nu (\ponto )
                           + t p (\ponto ) r_t (\ponto )
\, ,
\eformu
where $\nu\in L^2_0(p)$, $\{ r_t \}\subset  L^2_0(p)$ and
$r_t \stackrel{ L^2_0(p) } \setag 0$ as $t\downarrow 0$.
We have,  
\formu{bgr9a}
\frac{M(p_t) - M(p)}{t}
& = &
\frac{\int_\xxx x p_t (x) \lambda (dx) - \int_\xxx x p (x) \lambda (dx)}{t}
\\ \nonumber 
& = &
<\nu (\ponto ) , (\ponto )>_p + <r_t (\ponto ) , (\ponto ) >_p
\\ \nonumber 
& \longrightarrow &
<\nu (\ponto ) , (\ponto )>_p
\, , \quad \mbox{as } t\downarrow 0
\, .
\eformu
The last convergence comes from the continuity of the inner product and
the $L^2 (p) $ convergence of the path remainder term to zero.
 
\noindent
Define the function
\formub
M^\bullet_p (\ponto ) 
=
(\ponto ) - \int_\xxx x p(x) \lambda (dx)
\, .
\eformu
Clearly, $M^\bullet_p$ is in $L^2_0(p)$ and
\formu{bgr9}
< \nu , M^\bullet_p >_p 
=
< \nu (\ponto ) , (\ponto ) - \int_\xxx x p(x) \lambda (dx) >_p
=
<\nu (\ponto ) , (\ponto ) >_p
\, .
\eformu
Since (\ref{bgr9}) and (\ref{bgr9a}) hold for any $L^2$
differentiable path, we conclude that $M$ is differentiable 
with respect to the
$L^2$ tangent set and $M^\bullet_p$ is a gradient of $M$. 
An argument based on subsequences (c.f. Labouriau , 1998) yields
the differentiability of $M$ with respect to the essential 
tangent set, \ie the mean functional is differentiable with is the strongest 
sense we can define in our setup.

Note that in this example \reff{bgr1} holds for any differentiable path
with tangent $\nu\in\ttt (p)$.
However, according to our definition of functional differentiability
it would be enough if the condition \reff{bgr1} holds for one path
with tangent $\nu$.
\eproof 
\end{example}


Let $\phi : \ppp^* \setag \re^q$ be a differentiable functional at
$p\in\ppp^*$ with gradient $\phi^\bullet_p :\xxx \setag \re$.
It follows immediately from the definition of gradient that a 
function  $\phi^\star_p :\xxx \setag \re$ in $L^2_0 (a)$ is also a 
gradient of $\phi$ at $p$ if and only if,
\formu{grr2}
\forall \nu\in \ttt (p) , \quad
<\nu, \phi^\bullet_p >_p = \, <\nu, \phi^\star_p >_p
\, .
\eformu
We conclude from the remark above that if $\phi^\bullet_p$ is a
gradient of $\phi$ at $p$ and $\xi\in \{ \ttt (p) \}^\perp$
(\ie $\xi$ is in the orthogonal complement of the tangent space
with respect to $L^2_0 (p)$), then $\phi^\bullet_p +\xi $ is also a
gradient of $\phi$ at $p$. Hence, in general the gradient of a differentiable 
functional is not unique.

A gradient $\phi^\bullet_p$ of a differentiable functional at $p\in\ppp^*$ 
is said to be a {\it canonical gradient} if 
$\phi^\bullet_p (\ponto ) \in \bar\ttt (p)$.
Here $\bar\ttt (p)$ denotes the $L^2$ closure of the space spanned by
$\ttt (p)$.
The following proposition shows that there exists only one canonical
gradient (apart from almost surely equal functions) and gives a recipe
to compute the canonical gradient, namely by orthogonal projecting  any
gradient onto $\bar\ttt (p)$.
We will see that the canonical gradient plays a crucial rule in the theory 
of estimation of functionals.

\begin{prop}
\label{grp1}
Let $\phi :\ppp^*\setag\re$ be a differentiable functional at $p\in\ppp^*$.
If $\phi^\bullet_p :\xxx\setag\re^q$ is a gradient of $\phi$ at $p$, then the 
vector formed by the orthogonal projection of components of $\phi^\bullet_p$
onto $\bar\ttt (p)$, say
$$
( \prod\{ \phi^\bullet_{1p} \vert \bar\ttt (p) \} , \dots ,
   \prod\{ \phi^\bullet_{qp} \vert \bar\ttt (p) \} )^T 
\, ,
$$
is also a gradient of $\phi$ at $p$.
Furthermore, if $\phi^\ast_p$ is another gradient of $\phi$ at $p$, then
$$
( \prod\{ \phi^\bullet_{1p} \vert \bar\ttt (p) \} , \dots ,
  \prod\{ \phi^\bullet_{qp} \vert \bar\ttt (p) \} )^T 
=
( \prod\{ \phi^\ast_{1p} \vert \bar\ttt (p) \} , \dots ,
  \prod\{ \phi^\ast_{qp} \vert \bar\ttt (p) \} )^T
\, ,
$$
$p$ almost surely.
\end{prop}
\proof
We prove the proposition for the case where $q=1$. The same argument
applied componentwisely
proves the case for $q\in N$, but with a more notation.
From the projection theorem we have the following orthogonal decomposition
\formub
\phi^\bullet_p = \prod\{ \phi^\bullet_p \vert \bar\ttt (p) \}
+ \prod\{ \phi^\bullet_p \vert \bar\ttt (p) \}
\, .
\eformu
Here $\bar\ttt^\perp(p)$ is the orthogonal complement of $\bar\ttt (p)$ in
$L^2_0 (p)$. Hence
\formub
\prod\{ \phi^\bullet_p \vert \bar\ttt (p) \}
=
\phi^\bullet_p - \prod\{ \phi^\bullet_p \vert \bar\ttt^\perp (p) \}
\, .
\eformu
Since $\prod\{ \phi^\bullet_p \vert \bar\ttt^\perp (p) \}$ is orthogonal to 
$\bar\ttt (p)$,
we conclude from (\ref{grr2}) that $\prod\{ \phi^\bullet_p \vert \bar\ttt (p) \}$
is a gradient.

Reasoning analogously we conclude that if $\phi^\ast$ is another gradient of
$\phi$ at $p$, then 
\formub
\prod\{ \phi^\ast_p \vert \bar\ttt (p) \}
=
\phi^\ast_p - \prod\{ \phi^\ast_p \vert \bar\ttt^\perp (p) \}
\, .
\eformu
is a gradient of $\phi$ at $p$.
From (\ref{grr2}), for all $\nu\in \bar\ttt (p)$
\formub
< \prod\{ \phi^\bullet_p \vert \bar\ttt (p) \} , \nu >_p
=
< \prod\{ \phi^\ast_p \vert \bar\ttt (p) \} , \nu >_p
\eformu
and hence, for all $\nu\in \bar\ttt (p)$,
\formu{gr3}
< \prod\{ \phi^\bullet_p \vert \bar\ttt (p) \} -  
\prod\{ \phi^\ast_p \vert \bar\ttt (p) \}
, \nu >_p = 0
\, .
\eformu

In particular (\ref{gr3}) holds for 
\formub
\nu = \prod\{ \phi^\bullet_p \vert \bar\ttt (p) \} -  
\prod\{ \phi^\ast_p \vert \bar\ttt (p) \}
\, ,
\eformu
which yields
\formub
  \| \prod\{ \phi^\bullet_p \vert \bar\ttt (p) \}  - 
   \prod\{ \phi^\ast_p \vert \bar\ttt (p) \} \|^2_{L^2(p)}
= 0
\, .
\eformu
We conclude that $\prod\{ \phi^\bullet_p \vert T(p,\ppp) \} =  
\prod\{ \phi^\ast_p \vert T(p,\ppp) \}$ $p$ almost surely.
\eproof

\begin{example}[Mean functional continued]
It can be shown that the tangent space of the model $\ppp$ given
by (\ref{grzuca}) at each $p\in\ppp^*$ is the whole space $L^2_0 (p)$.
Hence the gradient calculated in example \ref{exgr2} is the canonical
gradient. Moreover, the canonical gradient is the only possible gradient
for the mean functional.
Note that if we drop the condition that requires the existence of the 
variance of $p$ (\ie condition \reff{bgr6}), then $M^\bullet_p$  
is no longer a gradient (because it is not in $L^2$) and $M$ is
not differentiable at $p$.
\eproof 
\end{example}

We consider next a proposition given trivial (but useful) rules for calculating
gradients of ``composed'' gradients.
\begin{prop}
\label{pprop8}
Let $\Psi,\phi :\ppp\setag\re^q$ be two differentiable functionals with
(canonical) gradient at $p\in\ppp^*$ $\Psi^\bullet$ and $\phi^\bullet$ respectively.
Let $g:\re^q \setag\re^q$ be a differentiable function.
\begin{description}
\item{ i)}
For all $a,b\in\re^q$, $a\Psi +b\phi$ is differentiable at $p$ and its
(canonical) gradient is given by $a\Psi^\bullet +b\phi^\bullet$
($a\Psi^\star +b\phi^\star $).
\item{ii)}
$g\circ\phi$ is differentiable at $p$ functional with gradient
$\nabla g\{ \phi (p) \} \{ \phi^\bullet (\ponto ) \}^T$.
If $\phi^\bullet$ is the canonical gradient of $\phi$ then 
$\nabla g\{ \phi (p) \} \{ \phi^\bullet (\ponto ) \}^T$
is the canonical gradient of $g\circ\phi$.
\end{description}
\end{prop}
\proof

\noindent
$i)$ Straightforward.

\noindent
$ii)$
We give next the proof for the case where $q=1$. 
The general case is obtained in a similar way.
Take an arbitrary differentiable path $\{ p_t \}$ with tangent $\nu$.
Define $\xi (t)=\phi (p_t)$, we have
\formub
\frac{\phi(p_t)-\phi(p)}{t}
\setag 
<\nu , \phi^\bullet >_p = \xi^\prime (0)
\, .
\eformu
Now,
\formub
\frac{ (g\circ\phi ) (p_t) - (g\circ\phi ) (p)}{t}
 \longrightarrow 
(g\circ\xi)^\prime (0) & = & g\left ( \xi (0) \right ) <\nu,\phi^\bullet >_p
\\ \nonumber & = & <\nu , \phi (p)\phi^\bullet >_p .
\eformu
\eproof

 
\subsubsection{Asymptotic bounds for functional estimation}

We study next some results concerning the estimation of a differentiable 
statistical functional under repeated sampling.
These results will illustrate the importance of the canonical gradient and
will guide the choice of the notion of path differentiability and
tangent cone to be used.

We start by defining sequences of estimators for a given 
differentiable functional 
$\phi :\ppp^*\longrightarrow\re^q$ (with respect to some tangent cones
$\ttt (p)$) based on samples.  
A sequence of functions $\{\hat\phi_n \}_{n\in N} = \{\hat\phi_n \}$
such that for each $n\in N$, $\hat\phi_n :\xxx^n\longrightarrow\re^q$ is
$(\aaa^n , \bbb (\re^q ) )$- measurable is said to be an 
{\it estimating sequence \/}.
Next we introduce two notions of regularity of estimating sequences often 
found in the literature.
An estimating sequence $\{\hat\phi_n \}$ is said to be {\it weakly regular\/} (for estimating $\phi$, with respect to the choice of tangent cones made)
if for each $p\in\ppp^*$ and each $\nu \in \ttt (p)$ there exists a 
differentiable path
$\{ p_{n^{-1/2}} \}_{n\in N}$ converging to $p$ and with domain
$V=\{ n^{-1/2} : n\in N \}$, for which
\formub
\sqrt{n} \{ \phi (p_{n^{-1/2}} ) - \phi (p) \}
\longrightarrow
\int_\xxx \phi^\bullet (x , p) \nu (x) p(x) \lambda (dx)
\eformu
and there exists a probability distribution $L_{p\nu}$
(not depending on the path) such that
\formub
\lll_{p^n_{n^{-1/2}}} 
\left [ \sqrt{n} \{ \hat\phi_n (\ponto ) \phi (p) \} \right ]
\stackrel{\ddd}\longrightarrow L_{p\nu}
\, .
\eformu
If the distributions $L_{p\nu}$ above do not depend on the tangent
$\nu$, then we say that $\{\hat\phi_n \}_{n\in N}$ is 
{\it regular}.

An important class of estimating sequences are the asymptotic linear 
sequences defined next.
An estimating sequence $\{\hat\phi_n \}$ is said to be {\it asymptotic linear\/}
(for estimating $\phi$) if there exists a function $IC_\phi :\xxx\times\ppp^*\longrightarrow\re$
such that for each $p\in\ppp^*$, the function $IC_\phi (\ponto ; p):\xxx\longrightarrow\re$
is in $L^2_0(p)$ and for each $n\in N$ given a sample 
${\bf x \/} = (x_1 ,\dots ,x_n)$ of size $n$, $\hat\phi_n$ admits the following 
representation
\formu{gr100}
\hat\phi_n ({\bf x \/} )
=
\phi (p) + \frac{1}{n} \sum_{i=1}^n IC_\phi ( x_i;p) 
+ 
o_{p^n} \left ( n^{-1/2} \right )
\, .
\eformu
The function $IC_\phi$ is called the {\it influence function\/} of $\phi$.
The representation (\ref{gr100}) can be re-written as
\formub
\sqrt {n} \left \{\hat\phi_n - \phi (p) \right \} =
 \frac{1}{\sqrt{n}} \sum_{i=1}^n IC_\phi ( x_i;p) + 
o_{p^n} \left (1 \right )
\, .
\eformu
From the central limit theorem and the Slutsky theorem
\formub
\sqrt {n} \left \{\hat\phi_n - \phi (p) \right \}
\stackrel{\ddd}\longrightarrow
N\left [ \zeb , Cov_p \{ IC_\phi (\ponto ;p) \} \right ]
\, ,
\eformu
where 
\formu
 Cov_p \{ IC_\phi (\ponto ;p) \} = \int_\xxx  IC_\phi (x;p) IC_\phi^T (x;p)
p(x)\lambda (dx)
\, .
\eformu
\begin{theor}
Let $\{\hat\phi_n \}$ be an asymptotic linear estimating sequence with 
influence function $IC$.
Suppose that for each $p\in\ppp^*$ the tangent cone is given by
$\ttt (p) = \stackrel{w} T^0 (p,\ppp^* )$.
Then,  $\{\hat\phi_n \}$ is regular if and only if for all $p\in\ppp^*$,
$\phi $ is differentiable at $p$ (with respect to $\ttt (p)$)
and $IC (\ponto ;p)$ is a gradient of $\phi$ at $p$.
\end{theor}
\proof
See Pfanzagl (1990) for the case where $q=1$ or Bickel {\it et al.\/} (1995).
\eproof

The theorem above identifies (influence functions of) regular asymptotic
linear sequences of estimators for estimating the functional $\phi$ with
the gradients of $\phi$.
The covariance, 
$\int \phi^\bullet_p (x) \phi^\bullet_p (x)^T p(x)\lambda (dx)$,
of a gradient $\phi^\bullet_p$ of $\phi$ is the asymptotic covariance of the 
corresponding regular asymptotic linear estimating sequence (under $p$) with 
influence function $\phi^\bullet_p$.
On the other hand, since the components of the canonical gradient $\phi^*$ of $\phi$ are the 
orthogonal projection of the components of any gradient onto the tangent space,
we have for a given gradient $\phi^\bullet_p$ and for all $p\in\ppp^*$
\formub
\phi^\bullet_p (\ponto )
=
\phi^\star_p (\ponto ) + R(\ponto ;p)
\, ,
\eformu
for some $R(\ponto ;p)\in \{ T^\perp (p;\ppp) \}^q$.
A standard argument yields then that, for all $p\in\ppp^*$,
\formu{marca111}
\hspace{-11mm}
\int_\xxx \hspace{-1mm} \left\{\phi^\ast _p(x) \right\} 
          \left\{\phi^\ast _p(x) \right\}^T p(x)\lambda (dx) \le \hspace{-2mm}
\int_\xxx \hspace{-1mm} \left\{\phi^\bullet_p (x) \right\}
          \left\{\phi^\bullet_p (x) \right\}^T p(x)\lambda (dx)
\, ,
\eformu
with inequality in the sense of the L\"owner partial order of matrices.
That is, the covariance of the canonical gradient is a lower bound
for the asymptotic covariance of regular asymptotic linear estimating
sequences.
Moreover, only an asymptotic linear estimating sequence 
with influence curve equal to the canonical gradient 
achieves this bound.
We say that an asymptotic linear estimating sequence is {\it optimal\/}
if, for each $p\in\ppp^*$, its influence function is the canonical gradient of $\phi$.
The bound \reff{marca111} is sometimes called the semiparametric 
Cram\`er-Rao bound.

In spite of the elegance of this theory, some care should be observed in 
applying it.
Firstly, there is a certain degree of arbitrariness in choosing only
the class of regular asymptotic linear estimating sequences.
When restricting to that class one can discard many interesting sequences.
This criticism applies, of course, to any optimality approach.
A second, more specific criticism is the following:
It occurs very often that the tangent space of large (semi- or non-parametric
models) is the whole space $L^2_0$ (see the examples at the end of the section
on tangent spaces).
In those cases, due to the uniqueness of the canonical gradient, each differentiable functional possesses only one gradient.
We conclude from the previous discussion that then there is only one possible influence function and hence all regular asymptotic linear estimating 
sequences are asymptotically equivalent (as far as the asymptotic variance is
concerned). Therefore an optimality theory for regular asymptotic linear
estimators is meaningless for the models with tangent spaces equal to the
whole $L^2_0$. 
  We refine next the optimality theory for functional estimation.

It is convenient to introduce the following notation.
Given a differentiable functional $\phi$ with respect to the tangent cones
$\{ \ttt (p) : p\in\ppp^* \}$ and with canonical gradient 
$\phi^\star (\ponto , p)$ at each $p\in\ppp^*$, denote 
$\int_\xxx \phi^\star _p(x) \phi^\star_p (x)^T p(x) \lambda (dx)$
by $I_\phi (p)$.
That is  $I_\phi (p)$ is the covariance matrix of the canonical gradient.
A weakly regular estimating sequence $\{ \hat\phi_n \}$ is 
{\it asymptotically of constant bias} at $p\in\ppp^*$ if for each 
$\nu , \eta \in \ttt (p)$
\formub
\int x d L_{p\nu} (x) 
=
\int x d L_{p\eta } (x)
\in\re^q
\, .
\eformu
In particular, any regular estimating sequence is asymptotically
of constant bias.
\begin{theor}
[van der Vaarts extended Cr\'amer-Rao theorem]
Let $\phi :\ppp^*\longrightarrow \re^q$ be a
differentiable at $p\in\ppp^*$ with respect to 
$\ttt (p)\subseteq \stackrel{w} T (p)$.
Suppose that the sequence $\{ \hat\phi_n \}$ is weakly regular
and asymptotically of constant bias at $p\in\ppp^*$.
Suppose also that the covariance matrix of $L_{p 0}$ exists.
Then
\formu{Cr1}
Cov (L_{p 0} ) \ge I_\phi (p)
\, ,
\eformu
where the symbol $''\ge ''$ is understood in the sense of the 
L\"owner partial order of matrices
\footnote{That is $A\ge B$ means that $A-B$ is positive definite.}.
Moreover, the equality in \reff{Cr1} occurs only if
\formu{Cr2}
\sqrt{n} \{ \hat\phi_n - \phi (p) \}
=
\frac{1}{\sqrt{n}} \sum_{j=1}^n \phi^\star_p (x_j ) + o_P (1)
\, .
\eformu
\end{theor}
\proof
See van der Vaart (1980).
\eproof
We see from the theorem above that the larger are the tangent cones
$\ttt (p)$ used, the sharper are the inequalities \reff{Cr1}.
Small tangent cones make more likely the differentiability of the
functional but can make also the bound in \reff{Cr1} unattainable.

Another important optimality result in the theory of estimation
of functionals is the convolution theorem, which we give the following
version.
\begin{theor}[Convolution theorem]
\label{theconv}
Suppose that $\ttt (p) $ is convex and $\phi :\ppp^*\longrightarrow \re^q$ 
differentiable at $p\in\ppp^*$ with respect to $\ttt (p)$. 
Then any limiting distribution $L_p$ of a regular estimating sequence
for $\phi$ at $p$ satisfies
\formu{conv1}
L_p = N ( \zeb ,I_\phi (p) ) * M
\, ,
\eformu
where $M$ is a probability measure on $\re^q$.
\end{theor}
\proof
See Pfanzagl (1990) for the case where $q=1$ and 
$\ttt (p)=\stackrel{w} T(p)$ and van der Vaart (1980) for the
general case.
\eproof
The expression (\ref{conv1}) shows that, under the assumptions of the 
convolution theorem, a regular estimating sequence cannot possess
asymptotic covariance smaller than the $L^2 (p)$ squared norm of the
canonical gradient.
This provides an extension of the interpretation of the optimality theory
for regular asymptotic linear estimating sequences.
In fact, even when the tangent cone is is the whole $L^2_0$, the ``optimal''
regular asymptotic linear estimating sequence attains the bound for the 
concentration of regular estimating sequences given by the convolution 
theorem, provide the functional is differentiable.
An advantage of the version of the convolution theorem presented is
that we need not to work with the whole tangent space but with a 
convex cone of it. This can be useful when the functional in study
is not differentiable or when the calculation of the (weak) tangent space
is not feasible.

We close this section presenting a theorem that gives a minimax 
approach to the problem of estimation of functionals.
A function $l:\re^q \longrightarrow \re$ is sad to be bowl-shaped if
$l(\zeb ) = 0$, $l(x)= l(-x)$ and for all $k\in\re$, $\{ x : l(x) \le k \}$
is convex.
\begin{theor}
[Local asymptotic minimax theorem]
\label{thelam1}
Suppose that for each $p\in\ppp^*$, $\ttt (p)\subseteq \stackrel{w} T (p)$is 
convex and $\phi :\ppp^*\longrightarrow \re^q$ differentiable
at $p\in\ppp^*$ with respect to $\ttt (p)$. Then

\noindent
i)
For any sequence of estimators which is weakly regular at $p$
and bowl-shaped loss function $l$
\formu{lam1}
\hspace{-8mm}
\sup_{\nu\in\ttt (p)}
\int l(x) dF_{p\nu} (x)
\ge
\int l(x) d N(\zeb , I_\phi (p) ) (x)
\, .
\eformu

\noindent
ii)
For any bowl-shaped loss function $l$ and any estimating sequence
$\{ \hat\phi_n \}$,
\formu{lam2}
\hspace{-9mm}
\lim_{c\rightarrow \infty} 
\liminf_{n\rightarrow\infty}
\hspace{-2mm}
\sup_{Q\in H_n (p,c) } 
\hspace{-2mm} 
E_Q \{ l [ \sqrt{n} \{ \hat\phi_n - \phi (Q) \} ] \}
\ge
\int l(x) d N (\zeb , I_\phi (p) ) \lambda (dx)
\, ,
\eformu
where $H_n (p , c) := \{ Q\in\ppp : n \int \{ dQ^{1/2} (x) - p^{1/2} (x)\}^2
                                                             \lambda (dx)\}$
is the interception between $\ppp$ and the ball constructed with the
Hellinger distance of center $p$ and radius $n^{-1/2}$.   
\end{theor} 
\proof
See van der Vaart (1980).
\eproof
Note that from part $i)$ one can obtain a bound for the concentration
of {\it weakly} regular estimating sequences based on the the canonical
gradient, provided $\phi$ is differentiable with respect to some convex
tangent cones.
In particular, if there exist an optimal asymptotic linear estimating
sequences and the assumptions of the theorem hold 
(\ie differentiability of $\phi$ and convexity of the tangent cone),
then the bound for weak regular estimating sequences given by \reff{lam1}
is attained by this regular asymptotic linear estimating sequence.
In this way,  in the case where the tangent space is the whole
$L^2_0$, the optimality of the (unique) regular asymptotic linear estimating
sequence can be justified.
The bound of the second part of the theorem above holds for the whole
class of estimators, however it is in general not attainable.
 
\newpage

\subsection {Asymptotic bounds for semiparametric models} 
\label{ggg4}

We consider a family of distributions $\ppp$ dominated by a 
the $\sigma$- finite measure $\lambda$ with representation
\formub
\ppp^*
=
\left  \{
\frac{d P_{\theta z}}{d \lambda } (\ponto )
=
p (\ponto ; \theta , z) 
\, : \,\,
\theta\in\Theta\subseteq\re^q 
\, , \,\,
z\in\zzz
\right \}
\, .
\eformu
Here $\theta$ is a $q$- dimensional interest parameter and 
$z$ is a nuisance parameter of arbitrary nature.
We assume that $\Theta$ is open and that the mapping 
$(\theta , z)\mapsto p(\ponto ;\theta , z)$
is a bijection between $\Theta\times\zzz$ and $\ppp^*$.
The {\it interest parameter functional } $\phi :\ppp^*\longrightarrow \re^q$ 
is defined, for each $p(\ponto ;\theta , z)\in\ppp^*$, by
\formub
\phi\{ p(\ponto ; \theta , z) \} = \theta 
\, .
\eformu
We will consider the differentiability of the interest parameter 
functional $\phi$ for a range of tangent cones.

Recall that we assumed that 
for each $(\theta_0 , z_0)\in\Theta\times\zzz$,
\formub
\forall x \in \xxx , \,\, p(x; \theta_0 , z_0 ) >0 \, ,
\eformu
that the partial score function
\formub
l (x;\theta_0 , z_0 ) = 
\frac{\nabla p(x ; \theta , z_0) \vert_{\theta = \theta_0}}
     {p(x ; \theta_0 , z_0)}
= \big (l_1 (x ; \theta_0 , z_0) , \dots , l_q (x ; \theta_0 , z_0)
  \big )^T
\eformu
is $\lambda$- almost everywhere well defined and
that for $i=1,\dots , q$,
\formub
l_i (x ; \theta_0 , z_0)\in L^2_0 (P_{\theta_0 z_0} )
\, .
\eformu

Let us consider a fixed $(\theta , z ) \in\Theta\times\zzz$ at which 
we will study the differentiability of $\phi$.
For notational simplicity we denote $p (\ponto ;\theta , z)$ by $p(\ponto )$.

The first tangent cone we consider is
\formub
\ttt_1 (p) = span \{ l_i (x ; \theta , z) : \, i=1,\dots ,q \}
\, .
\eformu
Take $\nu \in \ttt_1 (p)$.
There exists $\alphab \in \re^q$ such that 
$\nu (\ponto ) =  l^T (\ponto ;\theta , z) \alphab$.
Define (for $t$ small enough) the path
\formub
p_t (\ponto ) = p (\ponto ; \theta + t\alphab , z)
\, .
\eformu
Clearly, there exists $\{ r_t \}$ such that
\formu{semii1}
l^T (\ponto ;\theta , z)\alphab 
=
\frac {p (\ponto ; \theta + t\alphab , z) - p(\ponto )}{t p(\ponto )}
+ r_t (\ponto )
\, ,
\eformu
with $r_t (\ponto )\longrightarrow 0$ $\lambda$- almost everywhere.
Hence the path $\{ p_t \}$ is ($L^\infty$) differentiable with tangent
$l^T (\ponto ;\theta , z)\alphab$.
Moreover, 
\formub
\frac { \phi (p_t ) - \phi (p)}{t}
=
\frac {\theta + t\alphab - \theta}{t}
= 
\alphab
\, .
\eformu
Defining 
$\phi_p^\star (\ponto ) = Cov_{\theta z}^{-1} \{\, l(\ponto ;\theta , z)\,\}
l (\ponto ; \theta , z)$ we obtain,
\formub
\hspace{-10mm}
\int_\xxx  \phi_p^\star (x) \nu (x) p(x) \lambda (dx)
& = &
\int_\xxx 
\hspace{-2mm}
Cov_{\theta z}^{-1} \{ l(\ponto ;\theta , z) \}
l (x ; \theta , z) l^T (x ; \theta , z)\alphab p(x) \lambda (dx)
\\ \nonumber & =&
\alphab
= 
\lim_{t\rightarrow 0} \frac{\phi (p_t) - \phi (p)}{t}
\, .
\eformu
We conclude that $\phi$ is differentiable at $p$ with respect to
$\ttt_1 (p)$. Moreover,
\formub 
Cov_{\theta z} (\, l(\ponto ;\theta , z)\, )^{-1}
l (\ponto ; \theta , z)
\eformu
is the canonical gradient of $\phi$.
Note that we used (in \reff{semii1}) implicitly the $L^\infty$ path 
differentiability,
however the argument presented holds for any weaker path differentiability.
For, note that the essential point is that we identify 
(through \reff{semii1}) any element of the tangent cone $\ttt_1 (p)$
with a $L^\infty$ differentiable path. If we adopt a path differentiability
weaker than the $L^\infty$ differentiability, then the $L^\infty$ 
differentiable paths identified with the elements of the tangent cone
would be differentiable in the current sense also and the differentiability
of the functional $\phi$ follows from the argument presented above.

The efficient scores $I_\phi (p)$ (\ie the correlation matrix of the
canonical gradient of $\phi$ at $p$) is the inverse of the correlation 
matrix of the score function $l(\ponto ;\theta z)$.
The bounds for the asymptotic variance obtained with this naive choice
of tangent cones are not attainable in general. This will be apparent
from the development presented next where sharper bounds will be presented.

We introduce the notion of nuisance tangent space that plays a fundamental
rule in the estimation theory in semiparametric models.
For each $\theta_0\in\Theta$ consider the submodels
\formub
\ppp^*_{\theta_0} 
=
\{ p(\ponto ;\theta_0 , z ) \, : \,\, z\in\zzz \}
\, .
\eformu
The nuisance {\it tangent set} at $(\theta , z)\in \Theta\times\zzz$,
$T_N^0 (\theta , z)$,
is the tangent set of $\ppp^*_\theta$, \ie
$T_N^0 (\theta , z) = T^0 (p,\ppp^*_\theta )$.
The closure of the space spanned by the nuisance tangent set is called
the {\it nuisance tangent space} and denoted by $T_N (\theta , z)$.
Here we do not specify the notion of path differentiability adopted,
but when necessary a symbol will be superimposed. 

An alternative for the tangent cone better than $\ttt_1 (p)$ is
\formub
\ttt_2 (p) = span \{ l_i (x ; \theta , z) : \, i=1,\dots ,q \}
             \cup
             T_N^0 (\theta , z)
\, .
\eformu
We show next that $\phi$ is differentiable with respect to
$\ttt_2 (p)$, no matter which notion of path differentiability
we use.
Consider a $\nu\in T_N^0 (p) \subset \ttt_2 (p)$.
There is a differentiable path $\{ p_t \}$ contained in $\ppp_\theta^*$
with tangent $\nu$.
Since for each $t$, $p_t\in\ppp_\theta^*$, $\phi (p_t) = \theta = \phi(p)$
and
\formub
\frac{\phi (p_t) - \phi (p)}{t} = \zeb
\, .
\eformu
From the definition of functional differentiability, any gradient 
$\phi^\bullet_p$ of $\phi$ should satisfies, 
for each $\nu\in T_N^0 (\theta , z)$,
\formu{semi1}
\zeb = \lim_{t\searrow 0 } \frac{\phi (p_t) - \phi (p)}{t}
     = \int_\xxx \phi^\bullet_p (x) \nu (x) p(x) \lambda (dx)
\, .
\eformu
On the other hand, the argument presented in the case of the tangent
cone be $\ttt_1 (p)$ implies that, if 
$\nu\in span \{ l_i (x ; \theta , z) : \, i=1,\dots ,q \}$,  
say $\nu (\ponto ) = l(\ponto ; \theta , z)^T \alphab$, for some 
$\alphab\in\re^q$, then any gradient $\phi^\bullet_p$ of $\phi$
satisfies,
\formu{semi2}
\alphab = \int_\xxx \phi^\bullet_p (x) \nu (x) p(x) \lambda (dx)
\, .
\eformu
Clearly, the conditions \reff{semi1} and \reff{semi2} are sufficient to 
ensure that $\phi^\bullet_p$ is a gradient of $\phi$.
From these considerations, a natural candidate for being a gradient of
$\phi$ is the (standardised) projection of the score function onto the 
orthogonal complement of the nuisance tangent space.
Formally, define the function 
$l^E :\xxx\times\Theta\times\zzz\longrightarrow\re^q$ 
by, for each $(\theta , z )\in\Theta\times\zzz$,
   $l^E (\ponto ;\theta , z)=\left (l^E_1 (\ponto ;\theta , z), \dots ,
                                    l^E_q (\ponto ;\theta , z)\right )^T$
where, for $i=1,\dots , q$,
\formub
l^E_i (\ponto ;\theta , z)
= 
\prod (l_i (\ponto ;\theta , z) \vert T_N^\perp (\theta , z) )
\, .
\eformu
Here $\prod (g\vert A)$ is the orthogonal projection of 
$g\in L^2_0 (P_{\theta z})$ onto $A\subseteq L^2_0 (P_{\theta z})$.
Moreover, $T_N^\perp (\theta , z)$ is the orthogonal complement of
$T_N (\theta , z)$ in $L^2_0 (P_{\theta z})$.
The function $l^E$ is called the {\it efficient score function}
and we define the {\it efficient score} by
\formub
J(\theta , z) 
= 
\int_\xxx l^E(x ;\theta ,z) l^E(x ;\theta ,z)^T p(x) \lambda (dx)
\, .
\eformu

Define 
\formub
\phi^\star_p (\ponto ) 
=
J(\theta , z)^{-1} \, l^E(x ;\theta ,z)
\, .
\eformu
Clearly $\phi^\star_p $ satisfies \reff{semi1} and \reff{semi2}.
We conclude that $\phi^\star_p $ is a gradient of $\phi$. 
Moreover, $\phi^\star_p $ is the canonical gradient
(with respect to $\ttt_2 (p)$), since $\phi^\star_p $
is in the closure of the span of the tangent cone.

Note that choosing $\ttt_2 (p)$ as the tangent cone, the functional 
$\phi$ is still differentiable and we obtain a bound related with
the extended Cram\'er-Rao inequality sharper than the bound obtained
with $\ttt_1 (p)$.
However, since the $\ttt_2 (p)$ is not necessarily convex,
it is impossible to use the convolution theorem and the 
local minimax theorem.

A third alternative for the tangent cone is
\formub
\ttt_3 (p) & = & span \{ l_i (x ; \theta , z) : \, i=1,\dots ,q \}
                +
                T_N^0 (\theta , z)
\\ \nonumber
 & = &
\left  \{ l(\ponto ;\theta , z)^T \alphab  + \eta (\ponto ) \, : \,\,
       \alphab \in \re^q \, , \,\, \eta \in T_N^0 (\theta , z) 
\right \}
\, .
\eformu 
Clearly $\ttt_3 (p)$ is convex, however the functional $\phi$
is not necessarily differentiable.
We introduce next an additional assumption in the model that will
make $\phi$ differentiable.
Suppose that for each $\alphab\in\re^q$ and each 
$\eta\in T_N^0 (\theta , z)$ there exists a generalised sequence 
$\{ z_t \} = \{ z_t (\theta ,z) \}$ such that $\{ p_t \}\subset \ppp^*$,
given by
\formu{pppfan1}
p_t (\ponto ) = p (\ponto ; t\alphab + \theta , z_t )
\eformu
is a differentiable path with tangent 
$l^T (\ponto ;\theta ,z) \alphab + \eta (\ponto )$.
This assumption can be found often in the literature in an implicit
form (see for instance Pfanzagl, 1990, page 17, for the case where
$q=1$).
We prove differentiability of $\phi$ at $p$ with respect to $\ttt_3 (p)$
under \reff{pppfan1}.
Given 
$\nu (\ponto )= l(\ponto ;\theta , z)^T \alphab  + \eta (\ponto )\in\ttt_3 (p)$,
and taking a path $\{ p_t \}$ as in \reff{pppfan1} we obtain
\formub
\frac{\phi (p_t ) - \phi (p)}{t}
=
\frac{ t \alphab + \theta - \theta}{t} 
=
\alphab
\, .
\eformu
On the other hand,
\formub
\hspace{-11mm}
\int_\xxx \hspace{-3mm} J^{-1} (\theta , z) l^E \hspace{-1mm}(x;\theta , z) 
\nu (x) p(x) \lambda (dx)
& = &
\hspace{-3mm}
J^{-1} (\theta , z) 
\hspace{-2mm}
\int_\xxx \hspace{-3mm}l^E (x;\theta , z) l^T \hspace{-1mm}(x;\theta , z) p(x) \lambda (dx)
\alphab
\\ \nonumber & &
\hspace{-3mm}
+
J^{-1} (\theta , z) 
\hspace{-1mm} 
\int_\xxx \hspace{-1mm} l^E (x;\theta , z) \eta (x) p(x) \lambda (dx)
\\ \nonumber 
& = &
\alphab 
= 
\lim_{t\searrow 0} \frac{\phi (p_t ) - \phi (p)}{t}
\, .
\eformu
Hence $\phi$ is differentiable at $p$ with respect to $\ttt_3 (p)$ and 
\formub
\phi^\star_p (\ponto ) = J^{-1} (\theta , z)l^E (\ponto ;\theta , z)
\eformu
is the canonical gradient.
In other words, we obtained the same canonical gradient of $\phi$
if we work with $\ttt_2 (p )$ or $\ttt_3 (p)$ and
consequently the extended Cram\'er-Rao bound is also the same with the 
two choices of tangent cone.
Note that $\ttt_3 (p)$ is convex hence we can use the convolution and
the local asymptotic minimax theorems. 
This provides an additional justification of the extended Cram\'er-Rao bound
(via convolution theorem) and a optimality theory involving a larger class
of estimators, namely the weakly regular asymptotic linear estimating
sequences (as in the first part of the local asymptotic minimax theorem)
or even arbitrary estimating sequences (as in the second part of the local 
asymptotic minimax theorem).
However, we pay a price for these improvements, we have to introduce
regularity conditions on the model in order to obtain the differentiability
of the interest parameter functional.

It is current in the literature to take the whole (weak or Hellinger)
tangent set as the tangent cone, assume that the tangent set is
equal to $\ttt_3 (p)$ and use (implicitly) assumptions equivalent to
\reff{pppfan1} (see Pfanzagl, 1990 page 17). 
The strength of the approach based on tangent cones, and not necessarily on 
the whole tangent set, is that it allow us to graduate the regularity 
conditions.
We can avoid the assumptions mentioned above in the difficult
cases or take full advantage of them in the sufficiently regular cases.
The approach based on tangent cones allow us to treat the cases where the
tangent set is difficult (or virtually impossible) to calculate.

We conclude the section with a comment regarding reparametrisations.
Suppose that we reparametrise the model by considering the interest
parameter $g(\theta )$ instead of $\theta$.
Here $g$ is a one-to-one differentiable application from
$\re^q$ to $\re^q$.
The interest parameter functional becomes $g\circ \psi (P_{\theta z})
= g(\theta )$.
An application of the proposition \ref{pprop8} and the chain rule shows
that if an estimating sequence $\{ \hat\theta_n \}$ attains the 
semiparametric Cram\`er-Rao bound for estimating $\theta$
then the transformed sequence $\{ g(\hat\theta_n ) \}$ attains the 
Cram\`er-Rao bound for estimating $g(\theta )$.

\newpage


\section{Estimating and Quasi Inference Functions}
\label{cap3}
In this section the theory of inference functions for models with nuisance parameters is studied .
The basic definitions and properties of inference functions are given
in section \ref{estsect2}. There a related notion called quasi estimating 
function is also introduced. Quasi inference functions are essentially
functions of the observations, the interest parameter and (different
from the inference functions) of the nuisance parameter. They will provide
a way to formalise in a more clear way the theory of inference function
and relate inference functions with regular asymptotic linear estimators.
In order to construct an optimality theory for inference functions,
we define a class of what we call regular inference functions.
Two alternative (and equivalent) characterisations of the regular estimating 
functions are provided in the subsections \ref{sss22} and \ref{sss23}.
The second characterisation is motivated by differential geometric
considerations concerning the statistical model (inspired by Amari and
Kawanabe, 1996).

The characterisations referred to are used to derive an optimality theory
in section \ref{estsect3}. A necessary and sufficient condition for the
coincidence of the bound for the concentration of estimators based on estimating
functions and the semiparametric Cram\`er-Rao bound is provided in 
subsection \ref{sss5}.
This condition says essentially that the nuisance tangent space should
not depend on the nuisance parameter.

The last section contains some complementary material.
Subsection \ref{sss6} studies a technique for obtaining optimal
inference functions when   the likelihood
function can be decomposed in certain way.
In this way an alternative justification for the 
so called principle of conditioning will be provided.
A generalisation of the notion of inference function is introduced 
in subsection \ref{sss7}.
The section closes with a result that will allow us to characterise when
the semiparametric Cram\`er-Rao bound is attained by estimators derived
from regular inference functions. 

\subsection
         [Basic definitions and properties]
         {Estimating functions and quasi- inference functions: 
         basic definitions and properties}
\label{estsect2}
 
\subsubsection{Inference and quasi-inference functions}
\label{sss21}

A function $\Psi : \xxx\times\Theta \setag \re^q$ such that
for each $\theta\in\Theta$, the associated function  
$\Psi (\ponto ; \theta ,z ) : \xxx \setag \re^q$ is measurable,
is termed an {\it inference function\/}. 
Estimating functions are used to define sequences of estimators
for the parameter of interest $\theta$ in the following way.
Under a repeated independent sample squeme, given a sample
${\bf x} = (x_1 , \dots , x_n)^T$ of size $n$ of the (unknown) 
distribution $P_{\theta z} \in \ppp$, define $\widehat \theta_n$
implicitly by the solution of the equation
\formu{est1}
\sum_{i=1}^n \Psi (x_i ; \widehat\theta_n ) = \zeb
\, \, .
\eformu
Under regularity conditions each $\widehat\theta_n$ is well defined
and the sequence $\{ \widehat\theta_n \}$ is consistent (for estimating $\theta$)
and asymptotically normally distributed. We explore this fact to construct
an optimality theory.

We introduce next a notion related to inference functions.
A function $\Psi :\xxx\times\Theta\times\zzz\setag \re^q$, of the
parameters and the observations, such that for each $\theta\in\Theta$ and each
$z\in\zzz$, the function $\Psi (\ponto ;\theta , z) :\xxx\setag\re^q$ is
measurable is called a {\it quasi-inference function\/}.
Each inference function can be naturally identified with a quasi-inference function
by making it correspond to a suitable quasi-inference function constant on the nuisance
parameter. We make no distinction between inference functions and the 
corresponding quasi- inference functions. This abuse of language causes, in general,
no risk of ambiguity.

A quasi- inference function $\Psi :\xxx\times\Theta\times\zzz\setag\re^q$ 
such that the conditions (\ref{est11})-(\ref{est15a}) below are satisfied is said to 
be a {\it regular quasi-inference function \/}. The conditions  are,
with $\psi_i$ denoting the $i^{th}$ component of $\Psi$ and
for all $\theta_0\in\Theta$, all $z\in\zzz$ and all $i,j\in\{1,...,p\}$,
\formu{est11}
\psi_i (\ponto ; \theta_0, z ) \in L^2_0 (P_{\theta_0 z} );
\eformu
the partial derivative with respect to $\theta$ is well defined
(almost everywhere), \ie
\formu{est12}
\frac{\partial}{\partial \theta^j}\psi_i (\ponto ; \theta , z )
\left |_{\theta = \theta_0}\right . 
\mbox { exists };
\eformu
the order of differentiation with respect to $\theta$ and integration can be
exchanged in the following sense
\formu{est14}
 \hspace{-10mm}
\frac {\partial}
      {\partial \theta^j} \! \!
\int \! \psi_i (x;\theta , z) p(x;\theta , z) \lambda (dx) 
\! \left |_{\theta = \theta_0}\right .
\!\!
= 
\!
\int
\!
\frac {\partial}
      {\partial \theta^j} 
\left [
\psi_i (x;\theta , z) p(x;\theta , z) 
\right ]_{\theta = \theta_0}
\!\!
\lambda (dx) ;
\eformu
the following $q\times q$ matrix is nonsingular 
\formu{est15}
\mbox{E}_{\theta z}
\left\{ \nabla_\theta \Psi (\ponto ;\theta , z) \right\}
=
\left [ \int_\xxx \frac {\partial}{\partial \theta^j} \psi_i (x ; \theta , z)
\left |_{\theta = \theta_0} \right .
              p(x;\theta_0 z) \lambda (dx) \,\, \right ]_{i,j = 1, ... ,q} 
\hspace{-5mm};
\eformu
and
\formu{est15a}
\hspace{-10mm}
\mbox{E}_{\theta z}\!\!
\left\{ \!\Psi (\cdot ;\theta_0 , z)  \Psi^T\!\! (\cdot ;\theta_0 , z)\!\right\}
\! = \!
\left [ \int_\xxx \!\!\!\! \psi_i (x ; \theta_0 , z) \psi_j (x ; \theta_0 , z)
              p(x;\theta_0 z) \lambda (dx)\! \right ]_{i,j = 1, ... ,q} 
\eformu
is positive definite.

It is presupposed that the parametric partial score function is a regular
quasi-inference function.

 A regular quasi-inference function that does not depend on the 
nuisance parameter $z$ is said to be a {\it regular inference function}.

\subsubsection{First characterisation of regular inference functions}
\label{sss22}

In this section we give a characterisation of the class of 
regular inference functions.
\begin{prop}
\label{estT1}
Let $\Psi :{\cal X}\times \Theta \times\zzz\longrightarrow \re^q$
be a regular quasi- inference function with components $\psi _1,\dots ,\psi _q$.
For all $(\theta , z) \in \Theta\times\zzz $ and $i\in \{1,\dots ,q\}$, 
\begin{eqnarray}
\nonumber
\psi_i (\,\cdot\, ;\theta ,z ) \in  
T_N^\perp (\theta , z ) \,\, .
\end{eqnarray}
\end{prop}
Here and in the rest of this text  $T_N(\theta ,z )=
\stackrel{2}T_N(\theta ,z )$ and
$T_N^{\perp }(\theta ,z )$ is the orthogonal
complement of the nuisance tangent space $\stackrel{2}T_N(\theta ,z )$ 
in $L^2_0 (P_{\theta z})$.

\proof
Take $(\theta , z) \in \Theta\times\zzz $ and $i\in \{ 1, \dots , k\}$ 
$z\in {\cal Z}$ fixed and $\nu \in T^0_N (\theta , z)$ arbitrary. 
We prove that $\nu $
and $\psi_i (\, \cdot \, ;\theta )$ are orthogonal in 
$L^2(P_{\theta z})$.
This implies the proposition, because of the continuity of
the inner product.

Let $\{ p_t \}_{t\in V}$ be a differentiable path at $(\theta , z)$ with tangent $\nu$
and remainder term  $\{ r_t \}_{t\in V}$.
Using (\ref{gendiff101}), for each $t\in V$, 
\begin{eqnarray}
\nonumber
\hspace{-10mm}
\langle \nu (\,\cdot\, ) , \psi_i ( \,\cdot\, ; \theta , z )
\rangle_{\theta z}
& = & \hspace{-3mm}
\langle [\{ p_t (\,\cdot\, ) - p (\,\cdot\, ;\theta ,z )\} / tp(\,\cdot\,
;\theta, z )
 ]
- r_t (\,\cdot\, ) , \psi_i (\,\cdot\, ;\theta ,z) \rangle_{\theta z}
\\ \nonumber & & \\ \nonumber 
& = & \hspace{-3mm}
\frac{1}{t}
\int_{\cal X} \! \psi_i (x ; \theta ,z) p_t (x) d \mu (x) -
\frac{1}{t}
\int_{\cal X} \! \psi_i (x ; \theta ,z) p (x;\theta,z ) d \mu (x)
\\ \nonumber & &
- \langle r_t (\,\cdot\, ) , \psi_i (\, \cdot \,  ; \theta ,z)
\rangle_{\theta z}
\\ \nonumber & & \\ \nonumber
& = & - \langle r_t (\,\cdot\, ) , \psi_i (\, \cdot \,  ; \theta ,z)
 \rangle_{\theta z}\, .
\end{eqnarray}
Since $r_t \stackrel{L^2 (P_{\theta z} )}{\longrightarrow }0$, from the
continuity of the inner product, we conclude that 
$$
\langle \nu (\,\cdot\, ) , \psi_i ( \,\cdot\, ; \theta ,z)
\rangle_{\theta z} = 0 \, . 
$$
\eproof

If the quasi- inference function does not depend on the nuisance parameter
(\ie it corresponds to a genuine inference function), then we can obtain
a sharper result.
\begin{prop}
\label{estT2}
Let $\Psi :{\cal X}\times \Theta\longrightarrow \re^q$
be a regular inference function with components $\psi _1,\dots ,\psi _q$.
For all $\theta \in \Theta$ and $i\in \{1,\dots ,q\}$, 
\begin{eqnarray}
\nonumber
\psi_i (\,\cdot\, ;\theta ) \in \bigcap_{z\in {\cal Z}} 
T_N^\perp (\theta , z ) \,\, .
\end{eqnarray}
\end{prop}
In fact, the proposition above holds for the class of quasi- estimating
functions with expectation invariant with respect to the nuisance parameter.
\proof
Take $\theta\in\Theta$ and $i\in \{ 1, \dots , k\}$ fixed and arbitrary $%
\xi\in {\cal Z}$ and $\nu \in T^0_N (\theta , \xi)$. We prove that $\nu $
and $\psi_i (\, \cdot \, ;\theta )$ are orthogonal in  $L^2
(P_{\theta z})$.

Let $\{ p_t \}_{t\in V}$ be a differentiable path at $(\theta , z)$ with 
tangent $\nu$ and remainder term  $\{ r_t \}_{t\in V}$.
Using (\ref{gendiff101}), for each $t\in V$, 
\begin{eqnarray}
\nonumber
\hspace{-10mm}
\langle \nu (\,\cdot\, ) , \psi_i ( \,\cdot\, ; \theta , z )
\rangle_{\theta z}
=
- \langle r_t (\,\cdot\, ) , \psi_i (\, \cdot \,  ; \theta ,z)
 \rangle_{\theta z}\, .
\eformu
Since $r_t \stackrel{L^2 (P_{\theta z} )}{\longrightarrow }0$, from the
continuity of the inner product, we conclude that 
\formub
\langle \nu (\,\cdot\, ) , \psi_i ( \,\cdot\, ; \theta ,z)
\rangle_{\theta z} = 0  
\, .
\eformu
\eproof

\subsubsection{Amari-Kawanabe's geometric characterisation of regular 
            inference functions}
\label{sss23}

We present in this section a variant of the geometric theory of estimating
functions for semiparametric models given in Amari and Kawanabe (1996).
The development presented is closely connected with the theory given in that
paper, however it is not exactly the same. We point out the most
remarkable differences at the end of the section.

Take $(\theta , z)\in \Theta\times\zzz$ fixed. 
Given  $a \in L^2_0 (P_{\theta z} )$ and $z_* \in\zzz$ denote
$p(\ponto , \theta , z)$ and $p(\ponto , \theta , z_*)$ by
$p(\ponto ) $ and $p_* (\ponto )$ respectively and define the
{\it $m$-parallel transport} of $a$ from $z$ to $z_*$ by
\formub
\pam {z}{z_*} a (\ponto ) 
=
\frac{p(\ponto )}{p_* (\ponto )} a(\ponto )
\, .
\eformu
If $a$ posses a finite expectation under $P_{\theta z_*}$ we define the
{\it $e$-parallel transport} of $a$ from $z$ to $z_*$ by
\formub
\pae {z}{z_*} a (\ponto ) 
=
a(\ponto ) - \int_\xxx a(x) p(x;\theta , z_*) \lambda (dx)
\, .
\eformu
 
The basic properties of the $m$- and $e$-parallel transport are given next.
\begin{prop}
\label{thamari1}
We have for each $z,z_*\in\zzz$ and each $a,b\in L^2_0(p)\cap L^2_0(p_*)$:
\formu{paral1}
\int_\xxx \pam {z}{z_*} b (x) p_* (x) \lambda (dx)
=
\int_\xxx \pae {z}{z_*} b (x) p_* (x) \lambda (dx) 
=
0 
\, ;
\eformu
\formu{paral2}
\langle a , \pam {z}{z_*} b \rangle_{\theta z_*}
=
\langle a ,  b \rangle_{\theta z}
\, ;
\eformu
\formu{paral3}
\langle \pae {z}{z_*} a , \pam {z}{z_*} b \rangle_{\theta z_*}
=
\langle a ,  b \rangle_{\theta z}
\eformu
and
\formu{paral4}
\pam {z}{z_*} \pam {z_*}{z} a(\ponto ) 
=
\pae {z}{z_*} \pae {z_*}{z} a(\ponto )
=
\pam {z}{z} a(\ponto )
=
\pae {z}{z} a(\ponto )
=
a(\ponto )
\, .
\eformu
\end{prop}
\proof
Straightforward from the definition of $m$-parallel transport.
\eproof

The parallel transports defined above have their 
origin in differential geometric 
considerations for statistical parametric models ($\alpha$-connections).
We will not enter in details of the geometric theory for semiparametric
models, but refer instead to Amari and Kawanabe (1996) for an informal
discussion.
The parallel transports permit us to
change the inner product (see \reff{paral2}), \ie it permits us to
move from one $L^2$ space to another, keeping to certain extent the
structure given by the inner product of the first space.
For instance the $L^2$ orthogonality (\ie noncorrelation) is 
preserved after $m$-parallel transporting.
From the statistical viewpoint the $e$- and the $m$-parallel transport 
corresponds to correcting for the mean and correcting for the distribution,
respectively, when we move from one $L^2$ space to another.

The following class of functions will be of interest in the theory
of inference functions,
\formub
F_{IA} (\theta , z ) 
=
\left \{ 
\begin{array}{ll}
r\in T_N^\perp (\theta , z )
\, : \, &
\forall z_* \in\zzz \mbox{ and }
\forall \nu_*\in T_N (\theta , z_* ) , \,
\\ &
\langle \pam {z_*}{z} \nu_* , r \rangle_{L^2 (P_{\theta z})} = 0
\end{array}
\right \}
\, .
\eformu
When the $e$-parallel transport is well defined one can use
alternatively the relation
$\langle  \nu_* , \pae {z}{z_*}r \rangle_{L^2 (P_{\theta z_*})} = 0$
instead of 
$\langle \pam {z_*}{z} \nu_* , r \rangle_{L^2 (P_{\theta z})} = 0$.
Informally, $F_{IA}$ is the class of functions $r$ in $T_N^\perp (\theta ,z)$
such that $r$ corrected for the mean or corrected for the distribution
is orthogonal to each $T_N (\theta ,z_*)$ under $P_{\theta z_*}$
(for $z_*$ running in the whole $\zzz$).

\begin{prop}
\label{thamari2}
For each $(\theta , z)\in\Theta\times\zzz$,
$F_{IA} (\theta , z)$ is a closed subspace of $L^2_0 (P_{\theta z} )$.
\end{prop}
\proof
The linearity and the continuity of 
$\langle \pam {z_*}{z} \nu_* , ( \ponto ) \rangle_{L^2 (P_{\theta z})}$
implies that $F_{IA} (\theta , z)$ is a vector subspace and a closed
set in $L^2 (P_{\theta z_*})$, respectively.
\eproof

The following proposition gives a characterisation of regular estimating
functions in terms of the classes of functions $F_{IA}$'s.
\begin{prop}
\label{thamari3}
Given a regular inference function $\Psi$ with components
$\psi_1 , \dots ,\psi_q$, we have, for $i=1,\dots ,q$, for all
$\theta\in\Theta$ and all $z\in\zzz$,
\formub
\psi_i (\ponto , \theta ) \in F_{IA} (\theta , z)
\, .
\eformu
\end{prop}
\proof
Take $i=1,\dots ,q$, $\theta\in\Theta$ and all $z\in\zzz$ fixed.
Given any $z_*\in\zzz$ and $\nu_*\in T_N (\theta , z_* )$
we have from proposition \ref{estT2} that 
$\psi_i (\ponto ;\theta )\in T_N^{\perp_*} (\theta , z_* )$
and then 
\formub
\langle 
 \pam{z_*}{z} \nu_* , \psi_i (\ponto ;\theta )
\rangle_{L^2 (P_{\theta z})} = 0
\, .
\eformu
Since $z_*$ was chosen arbitrarily, 
$\psi_i (\ponto ;\theta )\in F_{IA} (\theta , z)$.
\eproof
The proposition above can be easily sharpened making the 
components of the regular inference functions belong to 
the intersection (over the nuisance parameter) of the $F_{IA}$'s.
However the following theorem shows that this is in fact not necessary,
since in fact $F_{IA}$ does not depend on the nuisance parameter.
We will use sometimes the notation $F_{IA} (\theta )$.
\begin{prop}
\label{thamari4}
For all $\theta\in\Theta$ and all $z\in\zzz$ we have,
\formub
F_{IA} (\theta , z) 
= 
\cap_{z_*\in \zzz} T_N^{\perp_*} (\theta , z_*)
\, .
\eformu 
\end{prop}
Here $T_N^{\perp_*} (\theta , z_*)$ denotes the orthogonal complement
of $T_N (\theta , z_*)$ in $L^2_0 (P_{\theta z_*})$.
\proof

\noindent
'$\subseteq$'
Take $\eta\in F_{IA} (\theta , z)$, $z_*\in\zzz$ arbitrary and 
$\nu_*\in T_N (\theta , z_* )$.
Applying \reff{paral2} yields
\formub
\langle \nu_* , \eta \rangle_{L^2 (P_{\theta z_*})}
=
\langle \pam{z_*}{z} \nu_* , \eta \rangle_{L^2 (P_{\theta z})} = 0
\, .
\eformu
Hence $\eta\in T_N^{\perp_*} (\theta , z_* )$.
Since $z_*$ was choose arbitrarily in $\zzz$,
$\eta \in \cap_{z_*\in \zzz} T_N^{\perp_*} (\theta , z_*)$.

\noindent
'$\supseteq$' 
Take an arbitrary $z_*\in\zzz$, 
$\eta\in\cap_{z_*\in \zzz} T_N^{\perp_*} (\theta , z_*)$
and $\nu_*\in T_N (\theta , z_*)$.
Using \reff{paral2} we obtain
\formub
\langle \pam{z_*}{z} \nu_* , \eta \rangle_{L^2 (P_{\theta z})}
=
\langle   \nu_* , \eta \rangle_{L^2 (P_{\theta z_*})} = 0
\, .
\eformu
Since $z_*$ is arbitrary, $\eta\in F_{IA} (\theta , z)$.
\eproof

The proposition \ref{thamari4} shows that the characterisation 
of regular inference functions obtained 
here is equivalent to what we obtained in the last section.
We remark that the characterisation based on the intersection
of the nuisance tangent spaces can be found in J\o rgensen and
Labouriau (1998) and the characterisation based on parallel transports
(\ie based on $F_{IA}$) is a variant of the results of Amari and
Kawanabe (1996). The main difference of the variant presented here and the 
original formulation in Amari and Kawanabe (1996) is that here we define 
via the $m$-parallel transport and there $F_{IA}$ is constructed through
$e$-parallel transport. Both formulations are equivalent from this
point of view, provided the $e$-parallel transport is well defined.
Moreover, when defining via the $m$-parallel transport the class $F_{IA}$
is automatically a closed subspace in $L^2_0$.

 
\subsection{Optimality theory for estimating func\-tions}
\label{estsect3}

\subsubsection{Classic optimality theory}

Given a regular (estimating) quasi-inference function $\Psi$ we define the 
{\it Godambe information} of $\Psi$ by
$J_\Psi :\Theta \times \zzz \longrightarrow \re^{2q}$,
where for each $\theta\in\Theta$ and each $z\in\zzz$,
\formu{est16}
\nonumber
\hspace{-8mm}
J_\Psi (\theta , z)\hspace{-1.0mm }
& = &\hspace{-0.9mm }
\mbox{E}_{\theta z} \!
\{ \nabla_\theta \Psi (\ponto ;\theta , z) \} 
\mbox{E}_{\theta z} \!\{\Psi (\ponto ; \theta , z ) 
\Psi^T \!(\ponto ;\theta , z) \}^{-1}
\mbox{E}_{\theta z}
\!
\{ \nabla_\theta \Psi (\ponto ;\theta , z) \}^T 
\\ \nonumber
 & = &
S_\Psi (\theta ,z) V^{-1}_\Psi (\theta z) S_\Psi^T (\theta ,z)
\, .
\eformu
Here 
\formub
S_\Psi (\theta ,z):=\mbox{E}_{\theta z}
\!
\{ \nabla_\theta \Psi (\ponto ;\theta , z) \}
\mbox{ and }
V_\Psi (\theta z):=
\mbox{E}_{\theta z} \!\{\Psi (\ponto ; \theta , z ) 
\Psi^T \!(\ponto ;\theta , z) \}
\eformu
are called the {\it sensibility} and the {\it variability}
of $\Psi$ (at $(\theta ,z)$ ), respectively.

Using  standard arguments based on a Taylor expansion of $\Psi$
it can be shown that under some additional regularity conditions
(each $\psi_i$ twice continuous differentiable with respect to
each component of $\theta$, for instance) a sequence 
$\{\widehat\theta_n \}$ of roots of a regular inference functions
is asymptotically normally distributed with asymptotic 
variance given by $J_\Psi^{-1} (\theta , z)$, provided 
$\{\widehat\theta_n \}$ is weakly consistent.
(see J{\o}rgensen and Labouriau, 1995 for conditions for consistency and
asymptotic normality).
Hence, we say that a
regular inference function $\Psi$ is {\it optimal\/}
when for all $\theta\in\Theta$, for all $z\in\zzz$ and for each
regular inference function $\Phi$,
\formub
J_\Phi (\theta , z ) \le J_\Psi (\theta , z ) \,\, .
\eformu
Here ''$\le$'' is understood in the sense of the L\"owner partial order of 
matrices given by the positive definiteness of the difference.

In the literature of inference functions it is customary to say that
it is possible to justify the use of some estimators using finite sample
arguments via inference functions and the Godambe estimation 
(see the articles of Godambe referred to).
The argument used there is that the Godambe information is a quantity
that should be maximised when using inference functions.
We do not share this point of view. The inference functions themselves
are not the object of our direct interest. Our concern with inference functions is
only through the estimators (or inferential procedures) associated with them.
Hence one should judge inference functions only through the properties of such
inferential procedures. In fact, apart from the asymptotic variance, there are no 
clear connections between the Godambe information and the 
(asymptotic or finite sample) properties of the estimators associated with
regular inference functions.

We say that two regular quasi-inference functions, 
$\Psi , \Phi : \xxx\times\Theta\times\zzz\setag\re^k$, 
are {\it equivalent} if, for each $\theta\in\Theta$ and $z\in\zzz$
there exists a $k\times k$ matrix with full rank $K(\theta , z )$,
such that
$$
\Psi (x ; \theta , z) = K(\theta , z) \Phi (x,\theta ,z ) 
\virg \,\,\,\, P_{\theta z} \mbox{ a.s. } \, .
$$
We stress that $K(\theta , z)$ must not depend on the observation $x$.
Clearly, two equivalent inference functions have the same roots almost
surely and hence produce essentially the same estimators, \ie
they are equivalent from the statistical point of view. 
Moreover, it is easy to see that two equivalent quasi-inference functions
share the same Godambe information for each value of the parameters.

\subsubsection{Lower bound for the asymptotic covariance of estimators obtained
            through inference functions}

We define the {\it information score function}, $l^I:{\cal X}\times \Theta
\times {\cal Z}\longrightarrow \re^q$, by the orthogonal projection of
the partial score function, $l$ onto $F_{IA} (\theta )$. More
precisely, for each $\theta \in \Theta $ and $z \in {\cal Z}$, the
$i$th component of the information score function ($i=1,\dots ,q$)
at $(\theta ,z )$ is given by 
\formub
l^I_i(\,\cdot\, ;\theta ,z )\,=\, \Pi \{l(\,\cdot\, ;\theta ,z )|F_{IA}
(\theta )\}, 
\, .
\eformu
 
The space spanned by the components $l_1^I,\dots ,l_q^I$ of the information
score function at $(\theta ,z )\in \Theta \times {\cal Z}$ is denoted by $%
E(\theta ,z )$, i.e. 
\formub
E(\theta ,z )=span\{l_i^I(\,\cdot \,;\theta ,z ):i=1,\dots ,q\}\,. 
\eformu
Note that $E(\theta ,z )$ is a closed (since it is finite-dimensional
vector space) subspace of $L_0^2(P_{\theta z })$. Hence given any regular
inference function $\Psi :{\cal X}\times \Theta \times %
\longrightarrow \re^q$ with components $\psi _1,\dots ,\psi _q$ we
have, for all $\theta \in \Theta $, $z \in {\cal Z}$ and $i\in \{1,\dots
,q\}$ the orthogonal decomposition 
\begin{eqnarray}
\label{F2004}
\psi_i (\,\cdot\, ; \theta ,z) =  \psi_i^A (\,\cdot\,  ; \theta ,z)
\, + \, \psi_i^I (\,\cdot\, ; \theta ,z)
\,\, ,
\end{eqnarray}
where $\psi_i^I(\,\cdot \,;\theta ,z )\in E(\theta ,z )$ and $\psi
_i^A(\,\cdot \,;\theta ,z )\in A(\theta ,z ):=E^{\perp }(\theta ,z )$.
Here $A(\theta ,z )$ is the orthogonal complement of $E(\theta ,z )$ in $%
L_0^2(P_{\theta z })$. The decomposition above induces the following
decomposition of each regular quasi-inference function 
\begin{eqnarray}
\label{F2005}
\Psi (\,\cdot\, ; \theta ,z) =  \Psi^A (\,\cdot\,  ; \theta ,z)
\, + \, \Psi^I (\,\cdot\, ; \theta ,z)
\,\, ,
\end{eqnarray}
where the components $\psi _i^A(\,\cdot \,;\theta ,z ),\dots ,\psi
_i^A(\,\cdot \,;\theta ,z )$ of $\Psi ^A$ at $(\theta ,z )$ are in $%
A(\theta ,z )$ and the components $\psi _i^I(\,\cdot \,;\theta ,z
),\dots ,\psi _i^I(\,\cdot \,;\theta ,z )$ of $\Psi ^I$ at $(\theta ,z )$
are in $E(\theta ,z )$.

We show next that taking the ``component'' $\Psi ^I$ of a regular (quasi-)
inference function improves the Godambe information. However, at this stage
a technical difficulty appears, the function $\Psi ^I$ is not necessarily a
regular quasi-inference function, and hence does not necessarily possesses a
well-defined Godambe information. For this reason we introduce next an
extension of the notion of sensitivity, and consequently of Godambe
information, which will make us able to speak of Godambe information of some
non-regular (quasi) inference functions. To motivate our extended notion of
sensitivity, consider a regular inference function $\Psi :{\cal X}\times
\Theta\times\zzz  \longrightarrow \re^q$. We characterise the
sensitivity of $\Psi $ in an alternative form that will suggest the
extension one should define. For each $(\theta ,z )\in \Theta \times {\cal %
Z}$ and each $i,j\in \{1,\dots ,q\}$ we have 
\begin{eqnarray}
0 & = & \frac {\partial}{\partial \theta_i} \int_{\cal X}
\psi_j (x;\theta ,z ) p(x;\theta ,z ) d\mu (x)
\\ \nonumber & & \mbox{ (differentiating under the integral sign )}
\\ \nonumber & = &
\int_{\cal X} \frac {\partial}{\partial \theta_i} \left \{
 \psi_j (x;\theta ,z ) p(x;\theta ,z )  \right \} d\lambda (x)
\\ \nonumber & = &
\int_{\cal X} \frac {\partial}{\partial \theta_i} \left \{
 \psi_j (x;\theta ,z ) \right \} p(x;\theta ,z )  d\lambda (x)
+
\int_{\cal X} \psi_j (x;\theta ,z )
\frac {\partial}{\partial \theta_i} \left \{p(x;\theta ,z )  \right \} d\lambda
(x)
\, .
\end{eqnarray}
Hence 
\begin{eqnarray}
\nonumber
 & & \int_{\cal X} \frac {\partial}{\partial \theta_i}
\left \{ \psi_j (x;\theta ,z )\right \} p(x;\theta ,z )   d\lambda (x)
\\ \nonumber & & =
- \int_{\cal X} \psi_j (x;\theta ,z )
\frac {\partial}{\partial \theta_i} \left \{p(x;\theta ,z )  \right \} d\lambda
(x)
\\ \nonumber &  & =
\int_{\cal X} \psi_j (x;\theta ,z )
l_i (x;\theta , z ) p(x;\theta ,z )  d\lambda (x)
=
- \langle \psi_j (\, \cdot \, ;\theta ,z ) , l_i (\, \cdot \,;\theta , z )
\rangle_{\theta z}
\\ \nonumber &  &
\mbox{(decomposing $l_i=l_i^A + l_i^I$ with $l_i^I\in F_{IA}$ and $l_i^A\in
T_N$)}
\\ \nonumber &  & =
- \langle \psi_j (\, \cdot \, ;\theta ,z ) , l_i^I (\, \cdot \,;\theta , z
)
\rangle_{\theta z}
- \langle \psi_j (\, \cdot \, ;\theta ,z ) , l_i^A (\, \cdot \, ;\theta , z
)
\rangle_{\theta z}
\\ \nonumber & &
\mbox{(Since $\psi_j\in F_{IA}$ and $l_i^A$ orthogonal $F_{IA}$)}
\\ \nonumber &  & =
- \langle \psi_j (\, \cdot \, ;\theta ,z ) , l_i^I (\, \cdot \, ;\theta , z
)
\rangle_{\theta z}
\\ \nonumber & &
\mbox{(decomposing  $\psi_j=\psi_i^A + \psi_i^I$ and using the orthogonality
       of $l_i^I$ and $\psi_i^A$)}
\\ \nonumber &  & =
 - \langle \psi_j^I (\, \cdot \, ;\theta ,z ) , l_i^I (\, \cdot \, ;\theta ,
 z )
\rangle_{\theta z}
\, .
\end{eqnarray}
We conclude that the sensitivity of $\Psi $ at $(\theta ,z )$ is given by 
\begin{eqnarray}
\label{F2010}
S_\Psi (\theta ,z )
=
\left [
 - \langle \psi_j^I (\, \cdot \, ;\theta ,z ) , l_i^I (\, \cdot \, ;\theta ,
 z )
\rangle_{\theta z}
\right ]_{i=1,\dots , k}^{j=1,\dots ,k}
\,\, .
\end{eqnarray}
Here $\left[ a_{ij}\right] _{i=1,\dots ,k}^{j=1,\dots ,k}$ denotes the
matrix formed by $a_{ij}$'s with $i$ indexing the columns and $j$ indexing
the lines.

We define the {\it extended sensitivity \/} (or simply the {\it sensitivity
\/}) of $\Psi $ by the matrix in the right-hand side of (\ref{F2010}). The
(extended) Godambe information is defined in the same way we did before but
using the extended sensitivity instead of the sensitivity. Note that both,
the standard and the extended, versions of the sensitivity (and the Godambe
information) coincide in the case where $\Psi$ is regular. Moreover, the
extended sensitivity is defined for each quasi-inference function whose
components are in $L^2_0$, not only for regular inference functions.
According to the new definition both $\Psi$ and $\Psi^I$ posses the same
sensitivity.

\begin{prop}
\label{T2010} Given a regular inference function $\Psi $, for all $\theta
\in \Theta $ and all $z \in {\cal Z}$, 
$$
J_\Psi (\theta ,z )\le J_{\Psi ^I}(\theta ,z )\,\,. 
$$
\end{prop}

\proof
For each $\theta\in\Theta$ and $z \in {\cal Z}$, 
\begin{eqnarray}
\nonumber
J^{-1}_\Psi (\theta ,z )
& = & S^{-1}_\Psi (\theta ,z ) V_\Psi (\theta ,z ) S^{-T}_\Psi (\theta ,z
)
\\ \nonumber & = &
S^{-1}_{\Psi^I} (\theta ,z )
\{ V_{\Psi^I} (\theta ,z ) + V_{\Psi^A} (\theta ,z ) \}S^{-T}_{\Psi^I}
(\theta ,z )
\\ \nonumber & = &
S^{-1}_{\Psi^I} (\theta ,z ) V_{\Psi^I} (\theta ,z )  S^{-T}_{\Psi^I}
(\theta ,z )
+
S^{-1}_{\Psi^I} (\theta ,z ) V_{\Psi^A} (\theta ,z )  S^{-T}_{\Psi^I}
(\theta ,z )
\\ \nonumber & \ge &
S^{-1}_{\Psi^I} (\theta ,z ) V_{\Psi^I} (\theta ,z )  S^{-T}_{\Psi^I}
(\theta ,z )
= J_{\Psi^I}^{-1} (\theta , z )
\, .
\end{eqnarray}
\eproof

The following proposition gives further properties of regular inference
functions, which will allow us to establish an upper bound for the Godambe
information.

\begin{prop}
\label{T2011} Given a regular inference function $\Psi $, for all $\theta
\in \Theta $ and all $z \in {\cal Z}$, we have:

\begin{itemize}
\item[(i)]  $\Psi ^I\sim l^I$;

\item[(ii)]  $span\{\Psi _i^I(\,\cdot \,;\theta ,z ):i=1,\dots
,k\}=E(\theta ,z )$;

\item[(iii)]  $J_{\Psi ^I}(\theta ,z )=J_{l^I}(\theta ,z )$.
\end{itemize}
\end{prop}

\proof
Take $\theta\in\Theta$ and $z\in {\cal Z}$ fixed.

\noindent
$(i)$ Assume without loss of generality that the components of the efficient
score function $l_1^I(\,\cdot \,;\theta ,z ),\dots ,l_q^I(\,\cdot
\,;\theta ,z )$ are orthonormal in $L_0^2(P_{\theta z })$. For each $%
i\in \{1,\dots ,q\}$, expanding $\psi _i(\,\cdot \,;\theta ,z )$ in a
Fourier series with respect to a basis whose first $q$ elements are $%
l_1^I(\,\cdot \,;\theta ,z ),\dots ,l_q^I(\,\cdot \,;\theta ,z )$ one
obtains 
\begin{eqnarray}
\nonumber
\psi_i (\,\cdot\, ;\theta ,z ) & = &
\langle l_1^I (\,\cdot\, ;\theta , z) , \psi_i (\,\cdot\, ;\theta ,z )
\rangle_{\theta z} \,
l_1^I (\,\cdot\, ;\theta , z)
\\ \nonumber & &
+ \cdots +
\langle l_q^I (\,\cdot\, ;\theta , z) , \psi_i (\,\cdot\, ;\theta ,z )
\rangle_{\theta z} \,
l_k^I (\,\cdot\, ;\theta , z)
+
\psi_i^A (\,\cdot\, ;\theta ,z )
\, .
\end{eqnarray}
That is, 
\begin{eqnarray}
\label{F2020}
\psi_i^I (\,\cdot\, ;\theta ,z ) & = &
\langle l_1^I (\,\cdot\, ;\theta , z) , \psi_i (\,\cdot\, ;\theta ,z )
\rangle_{\theta z} \,
l_1^I (\,\cdot\, ;\theta , z) \\ \nonumber & &
+ \cdots +
\langle l_q^I (\,\cdot\, ;\theta , z) , \psi_i (\,\cdot\, ;\theta ,z )
\rangle_{\theta z} \,
l_q^I (\,\cdot\, ;\theta , z)
\, .
\end{eqnarray}
Moreover, for $j=1,\dots ,q$ 
\begin{eqnarray}
\label{F2021}
\langle l_j^I (\,\cdot\, ;\theta , z) , \psi_i (\,\cdot\, ;\theta ,z )
\rangle_{\theta z}
=
- \int_{\cal X} \left \{ \frac {\partial}{\partial \theta_j} \psi_i (x ; \theta
, z )
		\right \}  p (x ; \theta , z ) d \lambda (x)
\, .
\end{eqnarray}

We conclude from (\ref{F2020}) and (\ref{F2021}) that $\Psi^I (\,\cdot \, ;
\theta ,z ) = S_\Psi (\theta ,z ) l^I (\,\cdot\, ;\theta ,z)$, which
means that $\Psi^I$ and $l^I$ are equivalent.

\noindent
$(ii)$ From the previous discussion $span\{\Psi _i^I(\,\cdot \,;\theta ,z
):i=1,\dots ,q\}$ is the space spanned by $-S_\Psi (\theta ,z )l^I(\,\cdot
\,;\theta ,z )$ which is the span of $\{l_i^I(\,\cdot \,;\theta ,z
):i=1,\dots ,q\}$, since the sensitivity by assumption is of full rank.

\noindent
$(iii)$ Straightforward.

\eproof

A consequence of the two last proposition is that $J_{l^I}$ is an upper
bound for the Godambe information of regular quasi inference functions. This
upper bound is attained by any (if any exists) extended regular inference
functions with components in $E$. In particular if $l^I$ is a regular
(quasi-) inference function, then it is an optimal (quasi-) inference
function. 

\subsubsection{Attainability of the semiparametric Cram\'er-Rao bound}
\label{sss5}

We study in this section the attainability of the semiparametric Cram\'er-Rao 
bound through regular inference function. 
More precisely, we give a necessary and sufficient
condition for the coincidence of the semiparametric Cram\'er-Rao bound and
the bound given in the previous section for the asymptotic variance of
estimators derived from regular inference functions.

Let us consider the interest parameter functional $\Phi :\ppp^*\setag\Theta$
given by, for each $(\theta , z)\in\Theta\times\zzz$,
\formub
\Phi (p(\ponto ;\theta , z)) = \theta
\, .
\eformu
As shown in chapter \ref{cap2}, the functional $\Phi$ is differentiable
at each $p(\ponto ;\theta ,z):=p(\ponto)\in\ppp^*$, with respect to the 
tangent cone  
$
\ttt (p) = T_N (\theta , z) 
\cup 
span \{ l_i (\ponto ;\theta ,z) : i=1,\dots ,q \}.
$
Here we adopt the $L^2$ path differentiability, since this is the path
differentiability used to characterise the class of regular estimating 
functions. 
We stress that the theory of functional differentiability used here
is compatible with any notion of path differentiability stronger than
(or equal to) the weak path differentiability, in particular the
$L^2$ path differentiability is allowed.
Moreover, in the examples we have in mind ($L^2$- restricted models) 
the notions of weak and $L^2$ path differentiability coincide.
Take a fixed $p(\ponto ;\theta , z)=p(\ponto )$ in $\ppp^*$.
Consider the function 
\formub
\phi^\bullet_p (\ponto )
=
Cov_p (l^I)^{-1} l^I (\ponto ; \theta , z)
\, .
\eformu
The following lemma will allow us to connect the optimality theory 
for inference functions with the semiparametric Cram\'er-Rao lower bound.
\begin{lemma}
\label{thqqq1}
The function $\phi^\bullet_p (\ponto )$ is a gradient
of $\Phi$ at $p$.
\end{lemma}
\proof
A little of reflection reveals that it is enough to verify that
\formu{qqq1}
\forall \nu\in T_N^0 (\theta ,z), \,\,
\int_\xxx \phi^\bullet_p (x) \nu (x) p(x) \lambda (dx)
=
\zeb
\eformu
and
\formu{qqq2}
\int_\xxx \phi^\bullet_p (x) l^T(x;\theta , z) p(x) \lambda (dx)
=
I_q
\, ,
\eformu
where $I_q$ is the $q\times q$ identity matrix.

Take $\nu\in T_N^0 (\theta ,z)$. Since each component of
$u^I (\ponto ;\theta , z)$ is in 
$F_{IA}(\theta , z)\subseteq T_N^\perp (\theta ,z)$,
condition \reff{qqq1} holds.
On the other hand,
\formub
\hspace{-10mm}
\int_\xxx \phi^\bullet_p (x) l^T(x;\theta , z) p(x) \lambda (dx)
& = &
\int_\xxx Cov_p (l^I)^{-1} l^I (x; \theta , z)  l^T(x;\theta , z) 
p(x) \lambda (dx)
\\ \nonumber & = &
Cov_p (l^I)^{-1}
\int_\xxx  l^I (x; \theta , z)  l^I(x;\theta , z)^T p(x) \lambda (dx)
\\ \nonumber & = & 
I_q
\, ,
\eformu
that is the condition \reff{qqq2} holds.
We conclude that $\Phi^\bullet_p$ is a gradient of $\Phi$ at $p$
with respect to the tangent cone $\ttt (p)$.
\eproof

According to the lemma above $\Phi^\bullet_p$ is a gradient of $\Phi$
at $p$, but not necessarily the canonical gradient.
In fact the canonical gradient of the functional $\Phi$ at 
$p$ with respect to $\ttt (p)$ is 
\formub
\Phi^\star_p (\ponto ) 
=
J^{-1} (\theta , z)
l^E (\ponto ;\theta , z)
\, ,
\eformu
where $l^E (\ponto ;\theta , z)$ is the efficient score function
at $(\theta ,z)$ and $J^{-1} (\theta , z)$ is the covariance matrix
of $l^E (\ponto ;\theta , z)$ under $P_{\theta z}$ 
(see chapter \ref{cap2}).
The unicity of the canonical gradient implies that
$\Phi^\bullet_p$ is the canonical gradient if and only if
it is equal to $\Phi^\star_p$ and this occurs if and only
if $T^\perp_N (\theta ,z)=F_{IA} (\theta)$.
The covariance of $\Phi^\star_p (\ponto )$ (under $P_{\theta z}$),
that is $J^{-1} (\theta , z)$,
gives the semiparametric Cram\'er-Rao lower bound. 
On the other hand, the lower bound for the asymptotic covariance of
estimators obtained from regular inference functions is the 
covariance (under $P_{\theta z}$) of $\Phi^\bullet_p$.
We conclude that the following result holds.
\begin{theor}
\label{thqqq2}
The semiparametric Cram\'er-Rao lower bound coincides with
the bound for the asymptotic covariance of estimators defined
through regular inference functions at $(\theta ,z)\in\Theta\times\zzz$
if and only if
\formub
T_N^\perp (\theta ,z)
=
F_{IA}(\theta )
(= \cap_{z_*\zzz}T_N^\perp (\theta ,z_*) )
\, .
\eformu
\end{theor}
The theorem above implies that inference functions produce
efficient estimators only if the orthogonal complement of the
nuisance tangent space does not depend on the nuisance parameter.


\subsection{Further aspects}
\label{estsect4}

\subsubsection{Optimal inference functions via conditioning}
\label{sss6}

We present in this section some results which allow us to compute optimal 
inference functions in many practical situations. 
The results will be in accordance with the so called conditioning 
principle.
For the sake of simplicity we study here only models with a one-dimensional 
parameter of interest. 

We study the situation where we have a likelihood factorisation of the
following form. Suppose that there exists a statistic $T\,=\,t(X)$ such
that, for all $\theta \in \Theta $ all $z \in {\cal Z}$ and all $x\in 
{\cal X}$, 
\begin{equation}
\label{F24}p(x;\theta ,z ) = f_t(x;\theta )h\{t(x);\theta ,z \}. 
\end{equation}

\begin{theor}
\label{T9}Assume that there
exists a statistic $T$ such that one has the decomposition (\ref{F24}).
Moreover, suppose that the class $\{P_{\theta z }^t:z \in {\cal Z}\}$,
where $P_{\theta z }^t$ is the distribution of $T(x)$ under $P_{\theta z
}$ ({\em i.e.\/} $X\sim P_{\theta z }$), is complete. Then the regular
inference function given by 
\begin{eqnarray}
\Psi (x;\theta )=\frac \partial {\partial \theta }\log f_t(x;\theta
),\,\forall x\in {\cal X},\forall \theta \in \Theta   \label{F25}
\end{eqnarray}
is optimal. Moreover, if $\Phi $ is also an optimal inference function then $%
\Phi $ is equivalent to $\Psi $.
\end{theor}

The theorem above gives an alternative justification for the use of
conditional inference.   

The following technical (and trivial) lemma will be the kernel of the proofs
that follow. But first it is convenient to introduce the following notation.
Given a regular inference function $\Psi: {\cal X}\times\Theta
\longrightarrow {I\!\! R}^k$, we define 
$$
\tilde{\Psi}(x; \theta )\,=\,\frac{\Psi(x;\theta ,z)}{\mbox{E}_{ \theta
z}\{\Psi^{\prime}( \theta )\}}, 
$$
which is called the standardised version of $\Psi$. Here $\Psi
^{\prime}(\theta) = \nabla_\theta \Psi (\theta )$.
Along this section we denote the class of all regular inference functions
by $\ggg$.

\begin{lemma}
\label{L4} For each regular
inference function $\Psi $ and $\Phi :{\cal X}\times \Theta \longrightarrow {%
\re}$ and each $\theta \in \Theta $ and $z \in {\cal Z}$, the following
assertions hold:

\begin{itemize}
\item[(i)]  
$$
\frac{\mbox{E}_{\theta z }\{\Psi (\theta )l(\theta ;z )\}}{\mbox{E} %
_{\theta z }\{\Psi ^{\prime }(\theta )\}}\,=\,-1, 
$$
where $l(\theta ;z )$ is the partial score function at $(\theta ;z )$;

\item[(ii)]  
$$
\mbox{E}_{\theta z }\left\{ \tilde \Phi ^2(\theta )\right\} \,=\,\mbox{E} %
_{\theta z }\left[ \{\tilde \Phi (\theta )-\tilde \Psi (\theta
)\}^2\right] \,+\,2\mbox{E}_{\theta z }\left\{ \tilde \Phi (\theta )\tilde
\Psi (\theta )\right\} \,-\,\mbox{E}_{\theta z }\left\{ \tilde \Psi
^2(\theta )\right\} \,\,. 
$$
\end{itemize}
\end{lemma}

\proof  Since $\Psi $ is unbiased, one has 
$$
\int \Psi (x;z )p(x;\theta ,z )d\mu (x)\,=\,0. 
$$
Differentiating the expectation above with respect to $\theta $ and
interchanging the order of differentiation and integration, we obtain 
$$
\mbox{E}_{\theta z }\left\{ \frac \partial {\partial \theta }\Psi (\theta
)\right\} \,+\,\mbox{E}_{\theta z }\{\Psi (\theta )l(\theta ;z )\}\,=\,0 
$$
which is equivalent to the first part of the lemma. The second part is
straightforward. \eproof

The following lemma gives a useful tool for computing optimal inference
functions.

\begin{lemma}
\label{T8} Assume the previous regularity conditions. Consider two functions 
$A:\Theta \longrightarrow {I\!\!R}\backslash \{0\}$ and $R:{\cal X}\times
\Theta \times {\cal Z}\longrightarrow {I\!\!R}$. Suppose that, for each
regular inference function $\Phi $, one has, for each $\theta \in \Theta $
and $z \in {\cal Z}$, 
$$
\int R(x;\theta ,z )\Phi (x;\theta )p(x;\theta ,z )d\mu (x)=0\,\,. 
$$

If a regular inference function $\Psi $ can be written in the form, for all $%
\theta \in \Theta $, 
\begin{equation}
\label{F20}\Psi (x;\theta )\,=\,A(\theta )l(x;\theta ,z )\,+\,R(x;\theta
,z ),
\end{equation}
for $x$ $P_{\theta z }$- almost surely, ($\Psi $ does not depend on $z $ 
even though $%
l$ and $R$ do), then $\Psi $ is optimal. Furthermore, a regular inference
function $\Phi $ is optimal if and only if for all $(\theta ,z)\in\Theta\zzz$,
$$
\tilde \Phi (\theta )\,=\,\tilde \Psi (\theta ) 
\,\,\,\,\,\,  , \mbox{ for $x$ $P_{\theta z}$ almost surely,}
$$
provided that there exists a decomposition as (\ref{F20}) above.
\end{lemma}

\proof  Take an arbitrary $(\theta ,z )\in \Theta \times {\cal Z}$. Given $%
\Phi \in {\cal G}$ one has 
\begin{eqnarray}
\mbox{E}_{\theta z }\{\tilde \Phi (\theta )\tilde \Psi (\theta )\} &\,=\,&%
\mbox{E}_{\theta z }\left[ \frac{\Phi (\theta )A(\theta )l(\theta ,z
)\,+\,\Phi (\theta )R(\cdot ;\theta ,z )}{\mbox{E}_{\theta z }\{\Phi
^{\prime }(\theta )\}\mbox{E}_{\theta z }\{\Psi ^{\prime }(\theta )\}} %
\right]   \label{F21} \\
&\,=\,&\frac{A(\theta )}{\mbox{E}_{\theta z }\{\Psi ^{\prime }(\theta )\}} %
\frac{\mbox{E}_{\theta z }\{\Phi (\theta )l(\theta ,z )\}}{\mbox{E}_{\theta z }\{\Phi ^{\prime }(\theta )\}}  \nonumber \\
&\,=\,&-\frac{A(\theta )}{\mbox{E}_{\theta z }\{\Psi ^{\prime }(\theta )\}} %
.  \nonumber
\end{eqnarray}
Hence the value of $\mbox{E}_{\theta z }\{\tilde \Phi (\theta )\tilde \Psi
(\theta )\}$ does not depend on $\Phi $, in particular, 
$$
\mbox{E}_{\theta z }\{\tilde \Phi (\theta )\tilde \Psi (\theta )\}\,=\,%
\mbox{E}_{\theta z }\{\tilde \Psi ^2(\theta )\}\,>0. 
$$
On the other hand, from (ii) of Lemma \ref{L4}, one has 
\begin{eqnarray}
\mbox{E}_{\theta z }\{\tilde \Phi ^2(\theta )\}\, &=&\,\mbox{E}_\theta
[\{\tilde \Phi (\theta )-\tilde \Psi (\theta )\}^2]\,+\,2\mbox{E}_{\theta
z }\{\tilde \Phi (\theta )\tilde \Psi (\theta )\}\,-\,\mbox{E}_{\theta z
}\{\tilde \Psi ^2(\theta )\}  \label{F22} \\
&=&\,\mbox{E}_\theta [\{\tilde \Phi (\theta )-\tilde \Psi (\theta )\}^2]\,+\,%
\mbox{E}_{\theta z }\{\tilde \Psi ^2(\theta )\}  \nonumber \\
&\geq &\,\mbox{E}_{\theta z }\{\tilde \Psi ^2(\theta )\},  \nonumber
\end{eqnarray}
for each $\Phi \in {\cal G}$. Thus, $\forall \theta \in \Theta ,\forall z
\in {\cal Z},\forall \Phi \in {\cal G}$, 
\begin{equation}
\label{F23}J_\Phi (\theta ,z )\,=\,\frac 1{\mbox{E}_{\theta z }\{\tilde
\Phi ^2(\theta )\}}\,\leq \,\frac 1{\mbox{E}_{\theta z }\{\tilde \Psi
^2(\theta )\}}\,=\,J_\Psi (\theta ,z ). 
\end{equation}
We conclude that $\Psi $ is optimal. For the second part of the theorem,
note that one has equality in (\ref{F22}), and hence in (\ref{F23}), if and
only if $\forall \theta \in \Theta ,\forall z \in {\cal Z,}$ $\mbox{E} %
_{\theta z }[\{\tilde \Phi (\theta )-\tilde \Psi (\theta )\}^2]\,=\,0$.
That is, if a regular inference function $\Phi $ is optimal then $\tilde
\Phi (\cdot ;\theta )\,=\,\tilde \Psi (\cdot ;\theta )\,\,P_{\theta z
}\mbox{-a.s. },\,\forall \theta \in \Theta ,\forall z \in {\cal Z}$. 
\eproof

We can prove now the main theorem of this section.
\proof
(of theorem \ref{T9})
Take $\theta\in\Theta$ and $z\in{\cal Z}$ fixed. From (\ref{F24}), 
\begin{eqnarray}  \label{F1010}
l(x;\theta,z) = \frac{\partial}{\partial\theta}\log p(x;\theta,z) = 
\frac{\partial}{\partial\theta}\log f_{t}(x;\theta) \,+\,\frac{\partial}{\partial\theta}\log h(x;\theta,z) \, .
\end{eqnarray}
We apply Theorem \ref{T8} to prove that $\psi$ is a (``unique") optimal
inference function. More precisely, defining $A(\theta)=1$ and $R(x;\theta
,z) =-\partial \log h \{ t(x) ;\theta ,z \} / \partial \theta$, and
using (\ref{F1010}) we can write $\Psi$ in the form 
\begin{eqnarray}
\Psi (x ; \theta ) = \frac{\partial}{\partial \theta} \log f_t (;\theta ) =
A(\theta ) l(x;\theta ,z) + R(x ; \theta ,z ) \, .  \nonumber
\end{eqnarray}
According to lemma \ref{T8}, if $R$ is orthogonal to every regular
inference function, then $\Psi$ is optimal, moreover $\Psi$ is the unique
optimal inference function, apart from equivalent inference functions.

Take an arbitrary regular inference function $\phi$. We show that $\phi$ and 
$R$ are orthogonal. Note that for each $z\in{\cal Z}$, 
$$
0 = \int \phi (x;\theta ) p(x;\theta ,z ) d\mu (x) = \int \phi (x;\theta )
f_t(x;\theta ) h\{ t(x) ;\theta ,z ) d\mu (x) \, . 
$$
On the other hand $\mbox{E}_{\theta z} (\phi \vert T) = \int \phi (x;\theta
) f_t(x;\theta )d\mu (x)$, which is independent of $z$. We write $\mbox{E} %
_{\theta}(\phi \vert T)$ for $\mbox{E}_{\theta z}(\phi \vert T)$, and we
have $\mbox{E}_{\theta z}\{E_{\theta}(\phi \vert T)\}=0$. Since $T$ is
complete, $E_{\theta}(\phi \vert T)=0$, $P_{\theta z}$ almost surely. We
have then, 
$$
\mbox{E}_{\theta z}\{\phi (\theta ) R(\theta,z )\} = \mbox{E} %
_{\theta z}\{ R(\theta,z ) E_{\theta}(\phi \vert T) \} = 0 \,. 
$$
\eproof


\subsubsection{Generalised inference functions}
\label{sss7}
 
A quasi inference function is said to be an 
{\it generalised inference function} if it is equivalent
to an inference function. 
More precisely, a quasi inference function 
$\Psi :\xxx\times\Theta\times\zzz\setag\re^q$
is an generalised inference function if for each 
$(\theta ,z)\in\Theta\times\zzz$ there exists a non-singular
$q\times q$ matrix $A(\theta ,z)$ and a measurable function
$\Phi (\ponto , \theta ):\xxx\setag\re^q$ such that
\formub
\Psi (x ;\theta , z) = A(\theta ,z)\Phi (x , \theta )
\, ,
\eformu
for $x$ $\lambda$- almost everywhere.
If $\Phi $ is a regular inference function, then
$\Psi $ is said to be a regular generalised inference functions.

Generalised inference functions are used for estimating the interest
parameter in the following way.
Given a sample $\xxb = (x_1 ,\dots ,x_n )^T$ of size $n$
of a unknown probability measure of the model, define the estimator
$\hat\theta_n$ implicitly by the solution of the following equation
\formub
\zeb & = & \sum_{i=1}^n \Psi (x_i ; \hat\theta_n ; z)
\\ \nonumber & = &
A(\hat\theta_n , z) \sum_{i=1}^n \Phi (x_i ; \hat\theta_n )
\, ,
\eformu
which is equivalent to 
\formub
\sum_{i=1}^n \Phi (x_i ; \hat\theta_n ) = \zeb
\, .
\eformu
In other words, for each generalised inference function there is
an inference function that yields the same estimating sequence.
Generalised inference functions are just a tool that will
simplify some formalisations.
Examples of generalised inference functions are the efficient estimating
function of most of the regular semiparametric models.
The following result will be useful latter on.
\begin{prop}
\label{ajadobo}
If for each $(\theta ,z)\in\Theta\times\zzz$,
$F_{IA} (\theta )=T_N^\perp (\theta ,z)$, 
and the efficient score function is a generalised estimating
function, then the inference function equivalent to the efficient
score function attains the semiparametric Cram\`er-Rao bound.
\end{prop}
\proof
Take an arbitrary $(\theta , z)\in\Theta\times\zzz$.
Since $F_{IA} (\theta )=T_N^\perp (\theta ,z)$, the efficient score function
$l^E$ coincides with the information score function $l^I$ at $(\theta ,z)$.
Hence the extended sensibility of the efficient score function 
(at $(\theta ,z)$) is the $q\times q$ identity matrix and the Godambe
information of $l^E$ at $(\theta ,z)$ is
\formub
J_{l^E} (\theta , z) = Cov_{\theta ,z}^{-1} (l^E )
\, ,
\eformu
which is the semiparametric Cram\`er-Rao bound.
If $\Phi$ is an inference function equivalent to the efficient score 
function then its Godambe information is equal to the Godambe information
of the efficient score function, that is $\Phi$ attains the semiparametric 
Cram\`er-Rao bound at $(\theta , z)$.
The proof follows now from the fact that $(\theta ,z)$ was chosen
arbitrarily.
\eproof


\end{document}